# Second-Order Multi-Set Allocation Occupancy (MAO) Distributions under Pairwise-Intersection Constraints: Exact Laws, MAO Norms, Inequalities, and Limit Theory

Xing-gang Mao

Department of Neurosurgery, Xijing Hospital, the Fourth Military Medical University, Xi'an, No. 17 Changle West Road, Xi'an, Shaanxi Province, China;

*Corresponding authors. E-mail addresses: xgmao@fmmu.edu.cn (X.G. Mao);

ORCID: 0000-0001-8505-1904

Short title: Second-Order MAO Distribution

Conflict of Interest Disclosures: None reported.

Funding/Support: This work was supported by the National Natural Science Foundation of China (Grant No: 82473137)

**Abstract**

Our previous work established the multi-set allocation occupancy (MAO) theory, including the MAO distribution, norms, inequalities, and limit theory, under constraints on marginal set sizes alone. By additionally conditioning on all pairwise intersections while leaving triple and higher-order membership patterns random, the present work develops its second-order extension. This extension naturally identifies the earlier marginal-size-only framework as the first-order MAO theory.

Let $N$ be the population size, $T$ the number of labelled sets, and $M$ a feasible symmetric matrix whose diagonal and off-diagonal entries specify marginal sizes and pairwise intersections. Write $a = (a_B)_{B\subseteq[T]}$ for the atomic membership counts, where $a_B$ is the number of elements belonging to exactly the sets indexed by $B$, and let $\mathcal{A}(N, M)$ be the resulting feasible atomic region. With labelled multiplicity

$$w(a) = \frac{N!}{\prod_{B\subseteq[T]} a_B\,!},\ Z_{N,M} = \sum_{a\in\mathcal{A}(N,M)} w\,(a),$$

we define the second-order exact-$t$ and at-least-$t$ MAO norms by the direct atom formulas

$$\parallel t^r \parallel_T = \frac{\sum_{a\in\mathcal{A}(N,M)}(\sum_{|B|=t} a_B)_r\, w(a)}{(N)_r Z_{N,M}},\ \|[t,T]^r\|_T = \frac{\sum_{a\in\mathcal{A}(N,M)}(\sum_{|B|\geq t} a_B)_r\, w(a)}{(N)_r Z_{N,M}}.$$

Based on the norms, the occupancy moments can be calculated exactly. If $X_{=t}$ and $X_{\geq t}$ denote the numbers of elements with membership degree exactly $t$ and at least $t$, respectively, then for every $\nu \geq 1$,

$$\mathbb{E}_{N,M}[X_{=t}^{\nu}] = \sum_{i=1}^{\nu} S_{\nu,i} \parallel t^i \parallel_T,\ \mathbb{E}_{N,M}[X_{\geq t}^{\nu}] = \sum_{i=1}^{\nu} S_{\nu,i} \left\|[t,T]^i\right\|_T.$$

We further developed Poisson and normal limiting theory for regular feasible sequences of second-order constraints, which are verified by numerical approximation. In addition, we show that repeated-vector MAO inequalities, valid in the first-order theory through averaging over the full symmetric configuration class, need not persist under second-order constraints. The theory is verified by computational routes, especially the brute-force enumeration. The model yields an exact null distribution conditional on all marginal membership frequencies and pairwise overlaps. It also provides a systematic route to higher-order MAO theories by prescribing intersections up to any chosen order.

## 1. Introduction

The number of elements shared by multiple sets is a fundamental statistic in probability, combinatorics, network analysis, genomics, ecology, epidemiology, and multi-experiment data integration[1-3]. Given $T$ labelled subsets of a finite population, two related questions occur frequently:

$$X_{=t} = \left|\left\{u \in [N]: \sum_{i=1}^{T} 1\{u \in A_i\} = t\right\}\right|,$$

the number of elements that occur in exactly $t$ sets, and

$$X_{\geq t} = \left|\left\{u \in [N]: \sum_{i=1}^{T} 1\{u \in A_i\} \geq t\right\}\right| = \sum_{s=t}^{T} X_{=s},$$

the number of elements that occur in at least $t$ sets. Here, $|S|$ denotes the cardinality of a finite set $S$.

The second statistic is often the more directly interpretable quantity. For example, after five repeated experiments, one may be interested in elements observed in at least four experiments rather than elements observed in exactly four experiments[3, 4]. Analogous examples include genes detected in at least $t$ studies, species observed in at least $t$ habitats or sampling periods, users exposed to at least $t$ campaigns, and nodes participating in at least $t$ communities or hyperedges.

Classical hypergeometric theory describes overlap between two subsets under fixed set sizes. Generalized hypergeometric distributions extend this perspective to any amount of subsets and study the distribution of elements overlapping at an exact or threshold level, which is also termed as the first-order multi-set allocation occupancy (MAO) distribution[1, 2]. The MAO framework provides a finite-population theory for these statistics under fixed marginal set sizes,

$$|A_i| = m_i,\ i \in [T],$$

including exact moments, MAO norms, MAO inequality and asymptotic approximations[5, 6].

A central limitation of first-order conditioning is that it does not preserve pairwise overlap structure. Set systems with identical marginal sizes may have very different degrees of attraction, exclusion, clustering, or co-occurrence. In many scientific settings, the pairwise intersections are observed quantities that should be retained rather than randomized away. This motivates conditioning simultaneously on

$$|A_i| = M_{ii},\ \ |A_i \cap A_j| = M_{ij},\ i < j.$$

We call the resulting model the **second-order MAO distribution**.

Venn diagrams, binary contingency tables, and fibers are primarily representational devices [7]. They encode the pattern-specific regions of a multi-set overlap system, or describe the configurations compatible with specified marginal constraints, but they do not by themselves yield a probability law for overlap multiplicity or a statistical benchmark for deciding whether an observed degree of repetition is exceptional. This distinction is essential: the full membership-pattern representation has $2^T$ cells, and the resulting combinatorial growth has long obstructed a tractable distributional theory beyond the classical generalized-hypergeometric setting.

MAO theory supplies this missing inferential layer. Its central random variable is the occupancy multiplicity

$$K = \sum_{i=1}^{T} 1\{U \in A_i\},$$

the number of sets containing a typical element. Rather than treating a multi-set overlap merely as a collection of Venn-diagram regions, MAO models the distribution of $K$ and the associated exact and threshold counts $X_{=t}$ and $X_{\geq t}$. It thereby provides a rigorous framework for asking whether the numbers of unique, repeatedly observed, or highly shared elements are statistically meaningful relative to an explicitly defined null structure.

The first-order MAO theory establishes the decisive foundation for this program. Its MAO norm turns occupancy multiplicity into a normalized statistical quantity and provides the structural template for higher-order extensions. The second-order MAO norm then incorporates pairwise-intersection information without abandoning the same occupancy-centered target. Thus, the first-order theory identifies the general form of the normalization, supplies a tractable source of second-order formulas, and supports a systematic progression toward higher-order MAO models. Overall, contingency tables resolve the geometry of overlap, and fibers describe its constrained realizations; MAO provides the missing probability theory for overlap multiplicity.

The main contributions of this paper are as follows. First, we define the second-order MAO distribution as a uniform conditional distribution over labelled set systems with prescribed marginal sizes and pairwise intersections. Second, we derive an exact atom-weighted partition function and full probability mass functions (PMFs) for both $X_{=t}$ and $X_{\geq t}$. Third, we introduce exact-$t$ and at-least-$t$ second-order MAO norms and show that they yield arbitrary factorial and ordinary moments. Fourth, we validate the formulas by independent exact computational methods. Fifth, we obtained the limit theory of the second MAO distribution, which reached to Poisson or normal distributions in different certain conditions.

Finally, we clarify the relation of the model to classical hypergeometric distributions, contingency-table conditional distributions, fixed-margin fibers, and higher-order interaction analysis.

## 2. Second-order MAO distribution

Let

$$[N] = \{1, \dots, N\}$$

be a labelled finite population. Consider $T$ labelled subsets

$$A_1, \dots, A_T \subseteq [N].$$

Let

$$M = \left(M_{ij}\right)_{1 \le i,j \le T}$$

be a symmetric nonnegative integer matrix satisfying

$$M_{ii} = | A_i |,$$

and

$$M_{ij} = | A_i \cap A_j |, \ i \neq j.$$

Thus, $M$ specifies the first-order set sizes through its diagonal entries and the second-order pairwise overlaps through its off-diagonal entries.

Define the constrained configuration space

$$\Omega(N, M) = \left\{ (A_1, \dots, A_T): \begin{array}{c} A_i \subseteq [N], \\ | A_i | = M_{ii}, \ i \in [T], \\ | A_i \cap A_j | = M_{ij}, \ 1 \le i < j \le T \end{array} \right\}.$$

We assume that $M$ is feasible, namely,

$$\Omega(N, M) \neq \emptyset.$$

The constrained configuration space $\Omega(N, M)$ contains all labelled $T$-tuples of subsets satisfying the prescribed first- and second-order overlap constraints. We place the uniform probability measure on this space:

$$\mathbb{P}_{N,M}(A_1, \dots, A_T) = \frac{1}{| \Omega(N, M) |}, \ (A_1, \dots, A_T) \in \Omega(N, M).$$

This is the MAO null model determined by $(N, M)$. Since $M$ is a symmetric matrix encoding pairwise intersections, it is the second-order instance of the general MAO framework. The associated MAO distributions are the induced distributions of the exact and threshold occupancy counts $X_{=t}$ and $X_{\geq t}$, as defined in the following.

For each element $u \in [N]$, define its membership order

$$K(u) = \sum_{i=1}^{T} 1\{u \in A_i\}.$$

The exact-$t$ occupancy count and the at-least-$t$ occupancy count are

$$X_{=t} = \sum_{u=1}^{N} 1\{K(u) = t\},$$

and

$$X_{\geq t} = \sum_{u=1}^{N} 1\{K(u) \geq t\} = \sum_{s=t}^{T} X_{=s}.$$

The case $t = 0$ is included for notational completeness:

$$X_{\geq 0} = N.$$

The constraints imply three deterministic identities:

$$\sum_{t=0}^{T} X_{=t} = N, \qquad (1)$$

$$\sum_{t=0}^{T} t\, X_{=t} = \sum_{i=1}^{T} M_{ii}\,, \qquad (2)$$

and

$$\sum_{t=0}^{T} \binom{t}{2} X_{=t} = \sum_{1\leq i<j\leq T} M_{ij}\,. \qquad (3)$$

Equation (1) counts all population elements. Equation (2) counts all element-set incidences. Equation (3) counts all pair-intersection incidences. Thus, the MAO distribution randomizes the higher-order membership structure while preserving exactly the prescribed first- and second-order overlap information encoded by $M$.

More generally, the same notation is retained when the type of $M$ changes. A vector $M$ specifies first-order constraints, a symmetric matrix $M$ specifies second-order constraints, and a higher-order tensor $M$ may specify higher-order overlap constraints. The model remains indexed by $(N, M)$.

## 3. Membership atoms and the exact partition function

For every subset $B \subseteq [T]$, define the membership atom

$$C_B = \left( \bigcap_{i \in B} A_i \right) \cap \left( \bigcap_{j \notin B} A_j^c \right),$$

and let

$$a_B = | C_B |.$$

Each $C_B$ contains those elements that belong to exactly the sets indexed by $B$. The sets $(C_B)_{B \subseteq [T]}$ divide the population into $2^T$ disjoint groups, one for each possible membership pattern. Every element belongs to one and only one such group. Hence, adding the sizes of all atoms gives the total number $N$ of population elements:

$$\sum_{B \subseteq [T]} a_B = N. \qquad (4)$$

The marginal set-size constraints become

$$\sum_{\substack{B \subseteq [T] \\ i \in B}} a_B = M_{ii}, \ i \in [T], \qquad (5)$$

and the pairwise-intersection constraints become

$$\sum_{\substack{B \subseteq [T] \\ \{i,j\} \subseteq B}} a_B = M_{ij}, \ 1 \le i < j \le T. \qquad (6)$$

Let $\mathcal{A}(N, M)$ be the collection of all nonnegative integer **atom vectors** satisfying (4)-(6):

$$\mathcal{A}(N, M) = \left\{ (a_B)_{B \subseteq [T]} \in \mathbb{Z}_{\ge 0}^{2^T} : (4), (5), \text{and } (6) \text{ hold} \right\}.$$

It is useful to distinguish three levels of description. Each labelled element $u \in [N]$ has a binary membership vector $Y_u = (\mathbf{1}\{u \in A_1\}, \dots, \mathbf{1}\{u \in A_T\}) \in \{0,1\}^T$. The atom count $a_B$ is the number of elements whose membership vector corresponds exactly to $B \subseteq [T]$. Hence the atom vector $a = (a_B)_{B \subseteq [T]}$ is a $2^T$-dimensional count vector, not a single binary membership vector. The set $\mathcal{A}(N, M)$ consists of all such count vectors compatible with the prescribed total size, set sizes, and pairwise intersections.

For a fixed feasible **atom vector** $a$, the counts $a_B$ specify how many labelled elements must be assigned to each membership atom, but not which labelled elements occupy those atoms. The number of compatible labelled assignments is therefore the multinomial coefficient $N! / \prod_B a_B!$.

For a fixed $M$, each feasible atom vector

$$a = (a_B)_{B \subseteq [T]} \in \mathcal{A}(N, M)$$

specifies one possible *atom-count configuration*: it records how many of the $N$ labelled population elements have each membership pattern $B$. It does not, however, specify which particular labelled elements occupy the atoms.

For a fixed atom-count configuration $a$, the number of labelled set systems that realize it is

$$w(a) = \frac{N!}{\prod_{B \subseteq [T]} a_B\,!}. \qquad (7)$$

Indeed, this is the number of ways to allocate the $N$ labelled elements among the $2^T$ membership atoms so that atom $C_B$ receives exactly $a_B$ elements. Once this allocation is fixed, the labelled sets $A_1, \ldots, A_T$ are uniquely determined.

The feasible labelled set systems are therefore grouped according to their atom-count configurations. Summing the number $w(a)$ of realizations over every feasible configuration gives the total number of labelled set systems consistent with $M$:

$$Z_{N,M} := \sum_{a \in \mathcal{A}(N,M)} w\,(a) = \sum_{a \in \mathcal{A}(N,M)} \frac{N!}{\prod_{B \subseteq [T]} a_B\,!}. \qquad (8)$$

Hence,

$$Z_{N,M} = \mid \Omega(N, M) \mid. \qquad (9)$$

In short, $w(a)$ counts the labelled realizations of one feasible atom-count configuration $a$, whereas $Z_{N,M}$ counts all labelled set systems satisfying the prescribed constraints $M$.

Equip $\Omega(N, M)$ with the uniform distribution; that is, choose a labelled set system uniformly at random from $\Omega(N, M)$. The induced distribution of its atom-count vector on $\mathcal{A}(N, M)$ is:

$$\mathbb{P}_{N,M}(a) = \frac{w(a)}{Z_{N,M}} = \frac{\frac{N!}{\prod_{B \subseteq [T]} a_B\,!}}{Z_{N,M}},\ a \in \mathcal{A}(N, M). \qquad (10)$$

Indeed, all labelled set systems are equally likely, and exactly $w(a)$ of them induce the atom-count vector $a$. Thus Atom-count vectors with larger multinomial weight $w(a)$ occur with greater probability because they are realized by more labelled set systems.

**4. Exact-$t$ and at-least-$t$ occupancy distributions**

Because every atom $B$ has membership order $\mid B \mid$,

$$X_{=t}(a) = \sum_{\substack{B \subseteq [T] \\ |B|=t}} a_B . \qquad (11)$$

Consequently, the exact finite probability mass function (PMF) of $X_{=t}$ is

$$\mathbb{P}_{N,M}(X_{=t} = x) = \frac{\sum_{\substack{a \in \mathcal{A}(N,M) \\ X_{=t}(a)=x}} \frac{N!}{\prod_{B \subseteq [T]} a_B !}}{\sum_{a \in \mathcal{A}(N,M)} \frac{N!}{\prod_{B \subseteq [T]} a_B !}}. \qquad (12)$$

The at-least-$t$ occupancy count is

$$X_{\geq t}(a) = \sum_{\substack{B \subseteq [T] \\ |B| \geq t}} a_B = \sum_{s=t}^{T} X_{=s}(a). \qquad (13)$$

Its exact PMF is therefore

$$\mathbb{P}_{N,M}(X_{\geq t} = x) = \frac{\sum_{\substack{a \in \mathcal{A}(N,M) \\ X_{\geq t}(a)=x}} \frac{N!}{\prod_{B \subseteq [T]} a_B !}}{\sum_{a \in \mathcal{A}(N,M)} \frac{N!}{\prod_{B \subseteq [T]} a_B !}}. \qquad (14)$$

The two distributions have different scientific interpretations. $X_{=t}$ isolates a single membership layer. In contrast, $X_{\geq t}$ aggregates all high-order membership layers and directly measures repeated, persistent, or broadly shared occurrence. For a five-experiment study, $X_{\geq 4}$ is precisely the number of elements reproduced in at least four experiments.

The boundary cases are

$$X_{\geq 0} = N, \qquad (15)$$

and

$$X_{\geq T} = X_{=T}. \qquad (16)$$

**5. MAO norms and arbitrary-order moments**

Let

$$(x)_i = x(x-1)\cdots(x-i+1)$$

denote the falling factorial, with

$$(x)_0 = 1.$$

Throughout this section, $\mathbb{E}_{N,M}$ and $\mathrm{Var}_{N,M}$ denote expectation and variance under the uniform MAO model determined by $N$ and the prescribed overlap structure $M$. In the present setting, $M$ is a symmetric matrix and hence specifies a second-order MAO model.

We retain the first-order MAO convention: MAO norms are defined on the factorial-moment scale. The parameters $(N, M)$ are suppressed in the norm notation, since they are fixed by the underlying MAO model.

For $i \geq 1$, define the $i$-fold exact-$t$ MAO norm by

$$\| t^i \|_T = \frac{\sum_{a\in\mathcal{A}(N,M)} \left(\sum_{\substack{B\subseteq[T]\\|B|=t}} a_B\right)_i \frac{N!}{\prod_{B\subseteq[T]} a_B\,!}}{(N)_i \sum_{a\in\mathcal{A}(N,M)} \frac{N!}{\prod_{B\subseteq[T]} a_B\,!}}. \tag{17}$$

Consequently we have,

$$\mathbb{E}_{N,M}[(X_{=t})_i] = \| t^i \|_T. \tag{18}$$

For at-least-$t$ occupancy, define

$$\left\|[t,T]^i\right\|_T = \frac{\sum_{a\in\mathcal{A}(N,M)} \left(\sum_{\substack{B\subseteq[T]\\|B|\geq t}} a_B\right)_i \frac{N!}{\prod_{B\subseteq[T]} a_B\,!}}{(N)_i \sum_{a\in\mathcal{A}(N,M)} \frac{N!}{\prod_{B\subseteq[T]} a_B\,!}}. \tag{19}$$

Hence,

$$\mathbb{E}_{N,M}[(X_{\geq t})_i] = \left\|[t,T]^i\right\|_T. \tag{20}$$

Equations (18) and (20) are the factorial-moment identities for the exact-$t$ and at-least-$t$ MAO occupancy distributions.

For every integer $v \geq 1$, the ordinary moments follow from

$$x^v = \sum_{1\leq i\leq v} s_{v,i}\,(x)_i,$$

where $s_{v,i}$ is a Stirling number of the second kind. Therefore,

$$\mathbb{E}_{N,M}[X_{=t}^v] = \sum_{1\leq i\leq v} s_{v,i}\left\|t^i\right\|_T, \tag{23}$$

and

$$\mathbb{E}_{N,M}[X_{\geq t}^{v}] = \sum_{1\leq i\leq v} s_{v,i} \left\| [t,T]^{i} \right\|_{T}. \qquad (24)$$

In particular,

$$\mathbb{E}_{N,M}[X_{=t}] = \| t \|_{T},$$

and

$$\mathbb{E}_{N,M}[X_{\geq t}] = \|[t,T]\|_{T}.$$

The variances are

$$\mathrm{Var}_{N,M}(X_{=t}) = \| t^{2} \|_{T} + \| t \|_{T} - \| t \|_{T}^{2}, \qquad (25)$$

and

$$\mathrm{Var}_{N,M}(X_{\geq t}) = \|[t,T]^{2}\|_{T} + \|[t,T]\|_{T} - \|[t,T]\|_{T}^{2}. \qquad (26)$$

Thus, the MAO moment calculus has the same form across all overlap orders: factorial moments are MAO norms, and arbitrary ordinary moments are obtained through the same Stirling transform.

The symbol $M$ denotes the prescribed overlap structure. When $M$ is a vector of set sizes, the model is the first-order MAO model. When $M$ is a symmetric matrix whose diagonal entries specify set sizes and whose off-diagonal entries specify pairwise overlaps, the model is the second-order MAO model considered here. More generally, a higher-order MAO model may be indexed by an appropriate overlap tensor $M$. The notation

$$\mathbb{P}_{N,M},\ \mathbb{E}_{N,M},\ \mathrm{Var}_{N,M},\ \Omega(N,M),\ \mathcal{A}(N,M),\ Z_{N,M},$$

remains unchanged across these cases.

**6. Second-Order MAO Norms as a Continuation of the First-Order Formulation**

For comparison with the first-order MAO formulation, fix the ordered (distinct) distinguished elements $1,\dots,r$. For membership signatures $S_1,\dots,S_r \subseteq [T]$, let $M^{S_1,\dots,S_r}$ denote the residual second-order margin obtained by removing the contributions of $S_1,\dots,S_r$, and define

$$g_{N,M}^{(2)}(S_1,\dots,S_r) := G_{N-r,\,M^{S_1,\dots,S_r}}^{(2)},\ G_{N,M}^{(2)} := |\ \Omega(N,M)\ | = Z_{N,M}.$$

with $g_{N,M}^{(2)}(S_1,\dots,S_r) = 0$ whenever the residual margin is infeasible.

Thus $g_{N,M}^{(2)}(S_1,\dots,S_r)$ counts the configurations in $\Omega(N,M)$ for which the distinguished elements have prescribed signatures $S_1,\dots,S_r$. For $p_1,\dots,p_r \in \{0,\dots,T\}$, define the corresponding transversal sum by

$$G^{(2)}_{N,M}(p_1,\dots,p_r) := \sum_{\substack{S_1,\dots,S_r\subseteq[T]\\ |S_1|=p_1,\dots,|S_r|=p_r}} g^{(2)}_{N,M}(S_1,\dots,S_r).$$

More generally, for $P_1,\dots,P_r \subseteq \{0,\dots,T\}$, set

$$G^{(2)}_{N,M}(P_1,\dots,P_r) := \sum_{\substack{S_1,\dots,S_r\subseteq[T]\\ |S_\ell|\in P_\ell\ (1\le\ell\le r)}} g^{(2)}_{N,M}(S_1,\dots,S_r),$$

and normalize by

$$\| (p_1,\dots,p_r) \|^{(2)}_{T;N,M} := \frac{G^{(2)}_{N,M}(p_1,\dots,p_r)}{G^{(2)}_{N,M}},\quad \| (P_1,\dots,P_r) \|^{(2)}_{T;N,M} := \frac{G^{(2)}_{N,M}(P_1,\dots,P_r)}{G^{(2)}_{N,M}}.$$

Indeed, if $X_{=t}(a) := \sum_{|B|=t} a_B$, then label symmetry gives

$$(N)_r\, G^{(2)}_{N,M}(t,\dots,t) = \sum_{a\in\mathcal{A}(N,M)} (X_{=t}(a))_r \frac{N!}{\prod_{B\subseteq[T]} a_B\,!}.$$

Since

$$G^{(2)}_{N,M} = \sum_{a\in\mathcal{A}(N,M)} \frac{N!}{\prod_{B\subseteq[T]} a_B\,!},$$

it follows that

$$\| t^r \|^{(2)}_{T;N,M} = \frac{G^{(2)}_{N,M}(t,\dots,t)}{G^{(2)}_{N,M}} = \frac{\sum_{a\in\mathcal{A}(N,M)} (X_{=t}(a))_r \frac{N!}{\prod_{B\subseteq[T]} a_B\,!}}{(N)_r \sum_{a\in\mathcal{A}(N,M)} \frac{N!}{\prod_{B\subseteq[T]} a_B\,!}}.$$

Considering that $X_{=t}(a) = \sum_{\substack{B\subseteq[T]\\ |B|=t}} a_B$, we get the formulas (17).

Likewise, writing $X_{\ge t}(a) := \sum_{|B|\ge t} a_B$, one obtains

$$\| [t,T]^r \|^{(2)}_{T;N,M} = \frac{G^{(2)}_{N,M}(\{t,\dots,T\},\dots,\{t,\dots,T\})}{G^{(2)}_{N,M}} = \mathbb{E}_{N,M}\,[(X_{\ge t})_r].$$

Hence these are exactly the $g/G$ representations of the atom formulas in (17) and (19). No additional factor $(N)_r$ appears in the $g^{(2)}/G^{(2)}$ ratios, because fixing $(1,\dots,r)$ has already divided out the $(N)_r$ possible ordered choices of distinct elements.

## 7. Distributional approximations and asymptotic questions

Exact enumeration provides the finite-sample distribution, factorial moments, ordinary moments, cumulants, and tail probabilities of MAO occupancy counts. However, exact computation may become expensive when $T$ is large or when the feasible atomic region $\mathcal{A}(N, M)$ has high dimension. Distributional approximations are therefore important for large constrained MAO systems.

The MAO factorial-moment calculus gives direct criteria for Poisson and normal limits. Let $M^{(N)}$ be a feasible sequence of second-order overlap matrices. Throughout, $T$ and $t$ are fixed as $N \to \infty$, and all probabilities, expectations, variances, norms, and cumulants are understood under the MAO model determined by $(N, M^{(N)})$.

### 7.1 Poisson approximation in the rare-occupancy regime

A Poisson approximation is appropriate when an occupancy class is rare and its expected count remains bounded. Recall that

$$\mathbb{E}_{N,M^{(N)}}[(X_{=t})_r] = \| t^r \|_T, \qquad (7.1)$$

and

$$\mathbb{E}_{N,M^{(N)}}[(X_{\geq t})_r] = \| [t, T]^r \|_T. \qquad (7.2)$$

In particular,

$$\mathbb{E}_{N,M^{(N)}}[X_{=t}] = \| t \|_T, \qquad (7.3)$$

and

$$\mathbb{E}_{N,M^{(N)}}[X_{\geq t}] = \| [t, T] \|_T. \qquad (7.4)$$

**Theorem 7.1. Poisson limit for MAO occupancy counts**

Let $X$ denote either $X_{=t}$ or $X_{\geq t}$. Suppose that, for every fixed integer $r \geq 1$,

$$\mathbb{E}_{N,M^{(N)}}[(X)_r] \longrightarrow \lambda^r \qquad (7.5)$$

for some $\lambda \in [0, \infty)$. Then

$$X \Rightarrow \text{Poisson}(\lambda). \qquad (7.6)$$

For the exact-$t$ occupancy count, condition (10.5) is

$$\| t^r \|_T \longrightarrow \lambda^r, \; r \geq 1. \qquad (7.7)$$

For the at-least-$t$ occupancy count, condition (8.5) is

$$\| [t, T]^r \|_T \longrightarrow \lambda^r, \; r \geq 1. \qquad (7.8)$$

The condition for $r = 1$ gives convergence of the mean:

$$\| t \|_T \longrightarrow \lambda \qquad (7.9)$$

in the exact-$t$ case, or

$$\|[t, T]\|_T \longrightarrow \lambda \qquad (7.10)$$

in the at-least-$t$ case. The higher-order conditions require that the MAO factorial moments asymptotically agree with the factorial moments of a Poisson random variable. Thus, bounded expectation alone does not imply a Poisson limit under arbitrary second-order MAO constraints.

### 7.2 Normal approximation for proportional MAO sequences

The moment and variance identities in (25) and (26) hold for every feasible pair $(N, M)$. We now consider a particularly natural asymptotic regime: a feasible second-order constraint matrix is enlarged by a common integer factor.

Let $N_0 \geq 1$, and let $M^{(0)}$ be feasible for a population of size $N_0$. For each integer $q \geq 1$, define

$$N_q := qN_0,\ M^{(N_q)} := qM^{(0)}. \qquad (7.11)$$

Thus, the set sizes and all pairwise intersection sizes grow proportionally with the population size. In particular, the normalized first- and second-order constraints remain fixed along the sequence.

We assume that this normalized constraint vector is an interior point of the feasible second-order moment region. Informally, this excludes boundary cases in which one or more membership patterns are forced to have asymptotically zero frequency. We also assume that the occupancy statistic under consideration is not completely determined by the fixed set sizes and pairwise intersections. We impose two regularity conditions. First, the normalized set sizes and pairwise intersections remain in the relative interior of the feasible second-order moment region. Equivalently, no membership pattern is forced to have asymptotically zero frequency. Second, the occupancy count under consideration is not already determined by the prescribed set sizes and pairwise intersections. Equivalently, its residual conditional variance is positive. The formal independent and identically distributed (iid) conditional representation and the corresponding lattice condition are given in Appendix C.

**Theorem 7.2 (Normal limit under proportional scaling).**

Fix $T$ and $t$. Consider the proportional feasible sequence $\left(N_q, M^{(N_q)}\right)$ defined in (7.11), and suppose that the regularity conditions of Appendix C hold.

1. For exact-$t$ occupancy, if the remaining conditional variance is nonzero, then

$$\operatorname{Var}_{N_q,M^{(N_q)}}(X_{=t}) = N_q v_{=t} + o(N_q) \qquad (7.12)$$

for some constant $v_{=t} > 0$, and

$$\frac{X_{=t} - \mathbb{E}_{N_q,M^{(N_q)}}[X_{=t}]}{\sqrt{\operatorname{Var}_{N_q,M^{(N_q)}}(X_{=t})}} \Rightarrow N(0,1). \qquad (7.13)$$

2. For at-least-$t$ occupancy, if the remaining conditional variance is nonzero, then

$$\operatorname{Var}_{N_q,M^{(N_q)}}(X_{\geq t}) = N_q v_{\geq t} + o(N_q) \qquad (7.14)$$

for some constant $v_{\geq t} > 0$, and

$$\frac{X_{\geq t} - \mathbb{E}_{N_q,M^{(N_q)}}[X_{\geq t}]}{\sqrt{\operatorname{Var}_{N_q,M^{(N_q)}}(X_{\geq t})}} \Rightarrow N(0,1). \qquad (7.15)$$

**Proof.** Appendix C represents the uniform MAO law with constraints $M^{(N_q)}$ as an iid maximum-entropy membership-vector law conditioned on its total first- and second-order statistics. Under proportional scaling, the conditioning value is exactly attainable for every $q$, and the underlying iid law does not change with $q$. The conditional lattice central limit theorem of Appendix C therefore applies to the sums defining $X_{=t}$ and $X_{\geq t}$. It yields the variance expansions (7.12) and (7.14), together with the normal limits (7.13) and (7.15).

The constants $v_{=t}$ and $v_{\geq t}$ measure the variation in the relevant occupancy indicators that remains after the first-order set sizes and second-order intersections have been fixed. Their positivity is necessary: if an occupancy count is already determined by these constraints, then its variance is zero and no nondegenerate normal limit is possible.

For example, when $T = 2$, the four membership-atom counts are determined by $N$, $M_{11}$, $M_{22}$, and $M_{12}$. Hence all exact and cumulative occupancy counts are deterministic in that case.

**7.3 Numerical verification.**

We illustrate two distinct asymptotic regimes for the highest-order count $X_{=3}$ in second-order MAO models. Under proportional scaling, the universe size and all entries of the second-order constraint matrix are multiplied by a common scale factor. In this regime, both $\mathbb{E}[X_{=3}]$ and $\operatorname{Var}(X_{=3})$ grow with the scale. The exact distribution becomes increasingly well approximated by the normal distribution with the same mean and variance, consistently with Theorem~7.2 (Fig. 1). In contrast, the Poisson distribution with matching mean remains inaccurate, because the exact variance is substantially smaller than the mean and the MAO distribution is correspondingly more concentrated than a Poisson law.

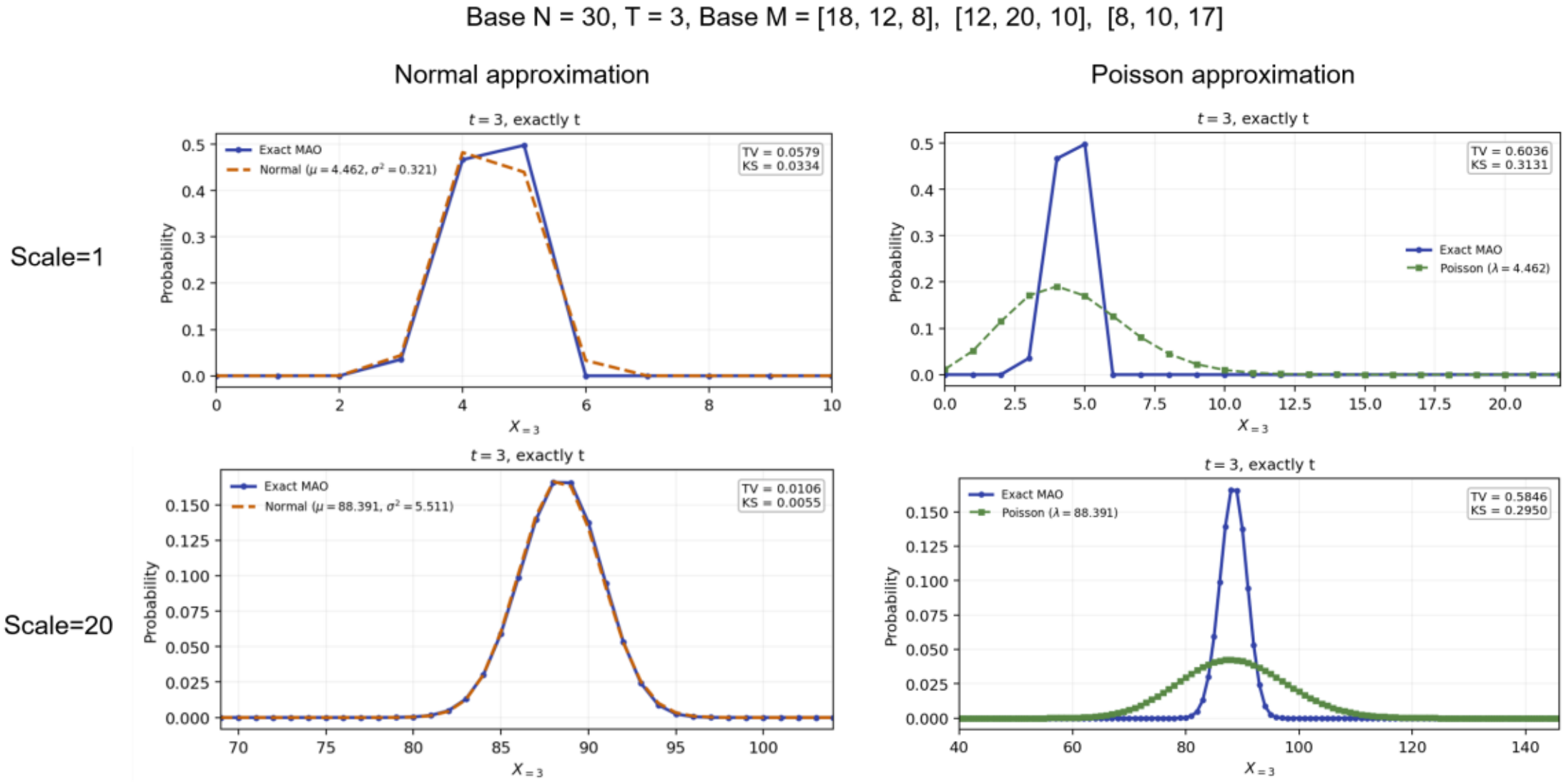

**Figure 1**. Normal and Poisson approximations under proportional scaling for a second-order MAO model with $T = 3$. The base instance has $N_0 = 30$, $M_0 = \begin{pmatrix} 18 & 12 & 8 \\ 12 & 20 & 10 \\ 8 & 10 & 17 \end{pmatrix}$. For each scale $q$, the universe size and every entry of the constraint matrix are scaled proportionally: $N^{(q)} = qN_0$, $M^{(q)} = qM_0$. The left column compares the exact distribution of $X_{=3}$ with the normal distribution having matching mean and variance, $\mathcal{N}\,(\mathbb{E}[X_{=3}], \mathrm{Var}(X_{=3}))$. The right column compares the exact distribution with the Poisson distribution having matching mean, $\mathrm{Pois}\ \ (\mathbb{E}[X_{=3}])$. At scale $q = 1$, $\mathbb{E}[X_{=3}] = 4.462$, $\mathrm{Var}(X_{=3}) = 0.321$. The normal approximation already has moderate accuracy, $\mathrm{TV} = 0.0579$, $\mathrm{KS} = 0.0334$, whereas the Poisson approximation is poor, $\mathrm{TV} = 0.6036$, $\mathrm{KS} = 0.3131$. At scale $q = 20$, $\mathbb{E}[X_{=3}] = 88.391$, $\mathrm{Var}(X_{=3}) = 5.511$. The exact distribution is then nearly indistinguishable from its moment-matched normal approximation, $\mathrm{TV} = 0.0106$, $\mathrm{KS} = 0.0055$. By contrast, the mean-matched Poisson approximation remains inaccurate, $\mathrm{TV} = 0.5846$, $\mathrm{KS} = 0.2950$. Thus, under proportional scaling, $X_{=3}$ becomes increasingly well approximated by a normal distribution, while its variance remains much smaller than its mean, $\mathrm{Var}(X_{=3}) \ll \mathbb{E}[X_{=3}]$, which rules out a Poisson approximation with matching mean.

We also consider a Poisson-type scaling in which the marginal constraints grow linearly with the scale, whereas the pairwise-overlap constraints are scaled so that $\mathbb{E}[X_{=3}]$ remains approximately constant. In the displayed example, the mean remains close to one. Although this mean is not vanishing in absolute magnitude, $X_{=3}$ remains sparse relative to the effective feasible range induced by the second-order constraints. As the scale increases, the dependence among the local contributions to $X_{=3}$ becomes weaker, and the mean-matched Poisson approximation improves substantially (Fig. 2). At the same time, the normal approximation does not exhibit a corresponding improvement.

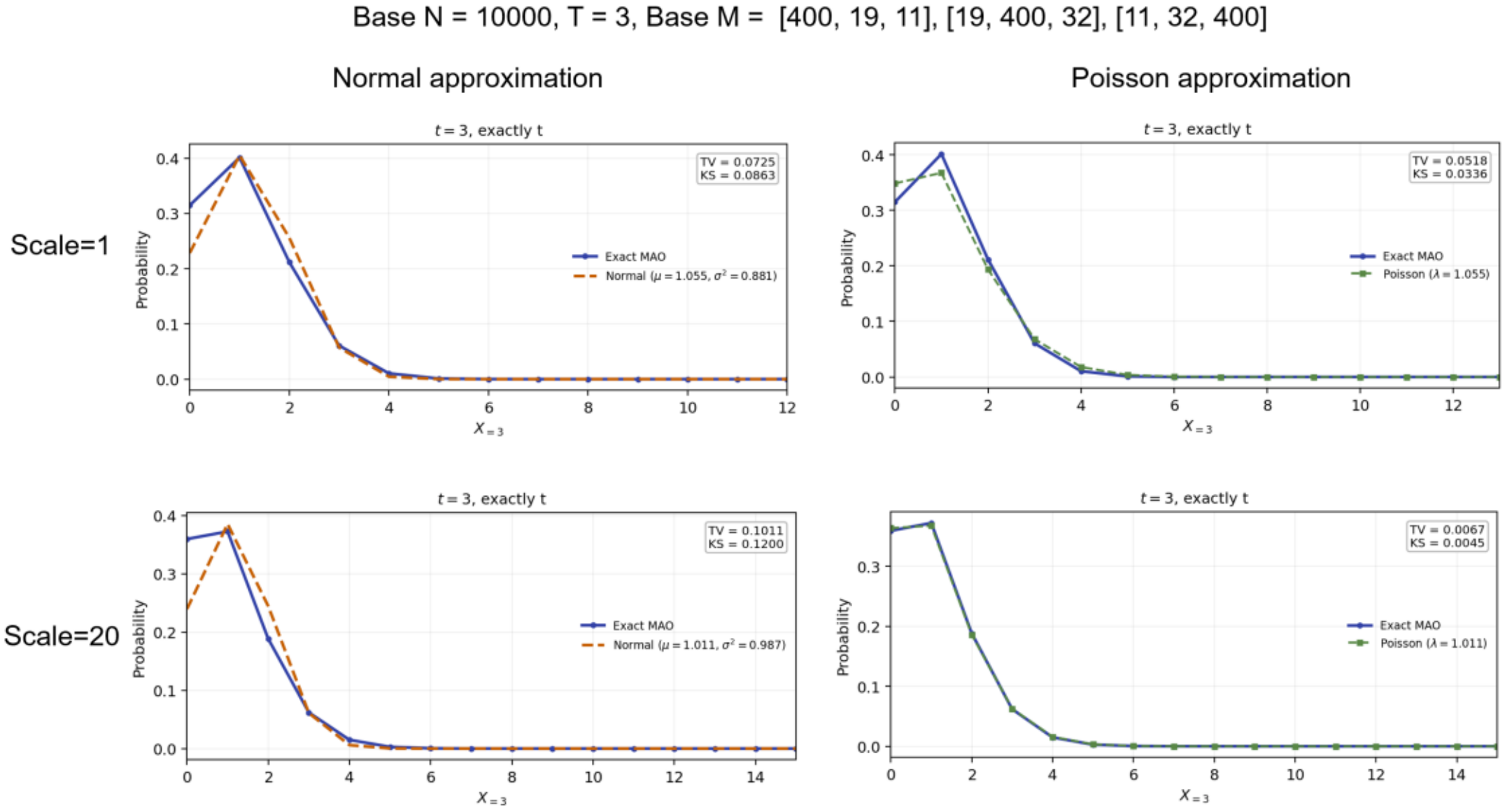

**Figure 2**. Normal and Poisson approximations under a Poisson-type scaling for a second-order MAO model with $T = 3$, based on $N_0 = 10{,}000$ and $M_0 = \begin{pmatrix} 400 & 19 & 11 \\ 19 & 400 & 32 \\ 11 & 32 & 400 \end{pmatrix}$. For each scale $q$, the universe size and marginal constraints scale linearly, $N^{(q)} = qN_0$ and $M_{ii}^{(q)} = q(M_0)_{ii}$, while pairwise-overlap constraints use $M_{ij}^{(q)} = \text{round}\left(q^{2/3}(M_0)_{ij}\right)$ for $i \neq j$. This scaling is designed to keep $\mathbb{E}\left[X_{=3}^{q}\right]$ approximately stable as $q$ increases; in the present example, it remains near one. The left column compares the exact MAO distribution with $\mathcal{N}\,(\mathbb{E}[X_{=3}], \text{Var}(X_{=3}))$, and the right column compares it with $\text{Pois}(\mathbb{E}[X_{=3}])$. At $q = 1$, $\mathbb{E}[X_{=3}] = 1.055$ and $\text{Var}(X_{=3}) = 0.881$; the normal approximation has $\text{TV} = 0.0725$ and $\text{KS} = 0.0863$, whereas the Poisson approximation has $\text{TV} = 0.0518$ and $\text{KS} = 0.0336$. At $q = 20$, $\mathbb{E}[X_{=3}] = 1.011$ and $\text{Var}(X_{=3}) = 0.987$; the normal approximation has $\text{TV} = 0.1011$ and $\text{KS} = 0.1200$, while the Poisson approximation improves to $\text{TV} = 0.0067$ and $\text{KS} = 0.0045$. Hence, as the scale increases, $\text{Var}(X_{=3})/\mathbb{E}[X_{=3}]$ approaches one and the exact distribution becomes nearly indistinguishable from the mean-matched Poisson law, illustrating the sparse-count Poisson regime induced by the second-order constraints.

These experiments provide numerical evidence for the two distributional regimes predicted by the theory: a normal regime when the highest-order count grows under proportional scaling, and a Poisson regime when its mean remains bounded while the relevant dependence becomes asymptotically weak. They further show that the distribution of $X_{=3}$ can be accurately estimated from the moments supplied by the second-order MAO norm calculations. In the normal regime, the mean and variance determine the moment-matched normal approximation; in the Poisson regime, the mean determines the Poisson approximation, while the comparison $\text{Var}(X_{=3}) \approx \mathbb{E}[X_{=3}]$ provides a diagnostic consistent with Poisson behavior. Thus, the second-order MAO framework links structural constraints, exact moment computation, and accurate distributional estimation.

**8. Repeated-vector inequalities under second-order constraints**

The repeated-vector MAO inequalities established in the first-order setting do not extend universally to the second-order model. In the first-order model, the feasible class consists of all set systems with prescribed set sizes. Equivalently, each set can be generated as the first prescribed number of positions of an independent uniform permutation of the $N$ labelled objects. The first-order inequality is consequently tied to averaging over this complete symmetric product class of permutations.

Prescribing pairwise intersections restricts this full first-order configuration class to a conditional subfamily with fixed overlap statistics. Although object-label exchangeability remains, this restricted family no longer retains the full combinatorial averaging structure underlying the first-order argument. Hence, the corresponding repeated-vector product bounds need not be preserved.

Exact rational computations confirm this distinction. While the inequalities hold in many second-order instances, explicit counterexamples occur for exact-occupancy norms. For example, in a feasible second-order instance,

$$\| 1,1 \|^{(2)} = \frac{1}{405} > \frac{1}{2025} = \left(\| 1 \|^{(2)}\right)^2.$$

Thus, the first-order inequality is not determined solely by the prescribed marginal and pairwise-overlap statistics. Its validity relies on the global combinatorial averaging over the full first-order configuration class, a structure that need not survive conditioning on pairwise intersections.

**9. Exact Computation and Validation**

For fixed $T$, exact computation is performed by enumerating the feasible atom vectors $a \in \mathcal{A}(N, M)$. The atom counts indexed by $B \subseteq [T]$ with $| B | \geq 3$ are treated as nonnegative integer variables subject to the residual pairwise-intersection constraints. Once these higher-order atoms have been assigned, the pair atoms are determined by the residual pairwise totals, the singleton atoms by the residual set-size constraints, and the empty atom by the population-size constraint.

Each feasible atom vector contributes the multinomial multiplicity

$$w(a) = \frac{N!}{\prod_{B \subseteq [T]} a_B\,!}$$

to the partition function

$$Z_{N,M} = \sum_{a \in \mathcal{A}(N,M)} w\,(a).$$

The same traversal accumulates the complete probability mass functions of the exact and threshold occupancy counts, $X_{=t}$ and $X_{\geq t}$, together with the factorial-moment numerators

$$\big(X_{=t}(a)\big)_r w(a), \qquad \big(X_{\geq t}(a)\big)_r w(a).$$

Thus, the occupancy PMFs, factorial moments, and second-order MAO norms are obtained from one exact finite sum. Partition functions, PMF weights, and factorial-moment numerators are accumulated using integer arithmetic. Probabilities, moments, and norms are formed only after normalization by $Z_{N,M}$.

Validation used three exact computational routes. First, atom-weighted enumeration produces the complete occupancy PMFs. Second, factorial moments and MAO norms are accumulated directly over the feasible atomic region. Third, for sufficiently small instances, exhaustive enumeration of labelled set systems in $\Omega(N, M)$ provides a combinatorially independent calculation. Agreement among these routes verifies the atom representation, multinomial weighting, occupancy distributions, factorial moments, and MAO norms.

As a representative example, consider $N = 9$, $T = 4$, and

$$M = \begin{pmatrix} 3 & 1 & 1 & 1 \\ 1 & 3 & 1 & 1 \\ 1 & 1 & 3 & 1 \\ 1 & 1 & 1 & 3 \end{pmatrix}.$$

Atomic enumeration yields

$$Z_{9,M} = 808920.$$

Direct enumeration of labelled set systems gives the identical count,

$$|\,\Omega(9, M)\,| = 808920.$$

As a second check, consider $X_{=3}$, the number of objects belonging to exactly three of the four sets. Direct enumeration gives

$$\Pr(X_{=3} = 0) = \frac{11}{107}, \ \Pr(X_{=3} = 1) = \frac{96}{107},$$

and hence

$$\mathbb{E}[X_{=3}] = \frac{96}{107}, \ \mathrm{Var}(X_{=3}) = \frac{1056}{11449}.$$

Direct atomic norm accumulation gives

$$\| 3 \|^{(2)} = \frac{32}{321}, \quad \| 3,3 \|^{(2)} = 0.$$

Consequently, the variance reconstructed from the first- and second-order MAO norms is

$$\mathrm{Var}(X_{=3}) = N \| 3 \|^{(2)} \left(1 - \| 3 \|^{(2)}\right) + N(N-1)\left(\|3,3\|^{(2)} - \| 3 \|^{(2)2}\right) = \frac{1056}{11449},$$

in exact agreement with the variance obtained from the directly enumerated PMF. The vanishing value $\| 3,3 \|^{(2)} = 0$ reflects the fact that no feasible set system in this benchmark can contain two distinct objects of exact membership order three.

Thus, the variance computed from the directly enumerated PMF agrees exactly with the variance reconstructed from the first- and second-order MAO norms. This validates, in particular, the multinomial weighting, the occupancy-count calculation, and the factorial-moment-to-norm correspondence.

Additional exact-validation results, including complete exact- and threshold-occupancy PMFs, factorial-moment reconstructions, second-order MAO norm comparisons, independent labelled-set enumeration checks, and higher-dimensional asymmetric benchmarks, are reported in Supplementary Material 2.

**10. Relation to existing distributions and conditional models**

The second-order MAO distribution is related to several established models, but it has a distinct inferential target and a distinct probability-space interpretation.

The classical hypergeometric distribution concerns the intersection of two subsets under fixed marginal sizes. It is recovered as the simplest overlap problem with $T = 2$. For $T > 2$, generalized hypergeometric distributions and first-order MAO distributions retain the finite-population occupancy perspective but condition only on the marginal set sizes[1, 2, 5]. They therefore average over all pairwise and higher-order overlap structures compatible with those sizes.

The second-order MAO distribution conditions further on every pairwise overlap. Thus, it is more restrictive than the first-order MAO model:

$$\text{first-order MAO: } \{| A_i |\}_{i=1}^{T} \text{ fixed,}$$

$$\text{second-order MAO: } \left\{|A_i| \ \left|A_i \cap A_j\right| \right\}_{1 \le i < j \le T} \text{ fixed.}$$

The second-order model can also be represented as a binary $N \times T$ incidence matrix. Each row describes the membership pattern of one labelled element, each column corresponds to a set, column sums equal $M_{ii}$, and column-pair inner products equal $M_{ij}$. In this representation, second-order MAO is a uniform distribution over labelled binary matrices with fixed column sums and fixed column-pair intersections.

At the level of atom counts, the model is related to a $2^T$-cell binary contingency table whose one-way and two-way margins are fixed. The atomic feasible region is a fixed-margin fiber. Such fibers are central objects in algebraic statistics, exact conditional inference, Markov-basis theory, and log-linear modeling [7-9].

However, the second-order MAO model is not simply a restatement of fixed-margin contingency-table theory. Its distinctive contribution is to provide a probability theory centered on **overlapping-set occupancy**. The objects of interest are not merely feasible contingency tables or tests of an interaction parameter. They are the exact distributions and arbitrary-order moments of interpretable occupancy statistics:

$$X_{=t},\ X_{\geq t},\ t = 0, \ldots, T.$$

This perspective has several advantages.

First, it directly targets the quantities reported in overlap studies: elements observed exactly $t$ times and elements replicated at least $t$ times. Second, it yields a unified factorial-moment representation through MAO norms. The same formalism covers all exact-$t$ and threshold occupancy statistics. Third, it links first-order and second-order conditioning within one nested framework. The first-order MAO distribution is recovered by averaging the second-order model over the pairwise-intersection matrices induced under first-order conditioning. Fourth, it separates fixed lower-order structure from random higher-order structure. When marginal sizes and pairwise overlaps are preserved, the distribution of $X_{\geq t}$ measures the amount of high-order recurrence or intersection compatible with those lower-order constraints. Finally, the model is naturally suited to exact finite-population null distributions. It does not require asymptotic independence, independent Bernoulli memberships, or unconstrained random graph assumptions.

## 11. Applications and potential impact

The second-order MAO distribution applies whenever multiple binary memberships are observed and both marginal frequencies and pairwise co-occurrences should be treated as fixed.

In repeated experiments, $X_{\geq t}$ quantifies reproducible signals. If an element is observed in at least $t$ of $T$ studies, the second-order MAO null distribution evaluates whether that degree of recurrence exceeds what is expected after preserving each experiment's yield and every pair of experiment overlaps[10, 11].

In genomics and multi-omics, sets may represent different gene lists, pathways, regulatory regions, cell populations, or experimental assays[12, 13]. Marginal sizes and pairwise overlaps are often strongly heterogeneous. The second-order MAO model provides a null distribution for high-order shared membership after retaining this heterogeneity.

In ecology, the sets can represent habitats, sites, seasons, species surveys, or environmental conditions[14, 15]. The count $X_{\geq t}$ measures the number of species present in at least $t$ locations or observation periods, conditional on observed richness and pairwise shared-species counts.

In epidemiology and clinical research, sets may encode diagnoses, exposures, medications, biomarkers, or eligibility criteria[16, 17]. A second-order MAO analysis can assess the expected number of individuals with multiple co-occurring conditions while preserving disease prevalences and pairwise comorbidity counts.

In network science and hypergraph analysis, each set may be a community, project, event, or hyperedge[18-21]. The set sizes are hyperedge degrees, the pairwise intersections are codegrees, and $X_{\geq t}$ counts nodes participating in at least $t$ groups. This is a direct high-order structural statistic under fixed degree and pairwise codegree sequences.

In database privacy and disclosure analysis, public release mechanisms may reveal marginal query counts and pairwise intersections[22]. The second-order MAO distribution quantifies the remaining uncertainty in higher-order overlaps and can be used to evaluate the probability of unusually large multi-query membership counts.

In recommendation, advertising, and behavior analysis, sets can represent audiences reached by campaigns or content categories[23]. The model distinguishes high-order exposure patterns attributable to observed pairwise campaign overlap from those indicating additional clustering beyond pairwise structure.

More generally, second-order MAO provides an exact microcanonical null model for high-order membership structure under prescribed first- and second-order sufficient statistics[20, 24, 25]. This places it at the intersection of finite-population probability, conditional inference, binary contingency-table analysis, and higher-order network science.

## 11. Conclusion

We developed the second-order multi-set allocation occupancy distribution for finite populations with fixed set sizes and fixed pairwise intersections. The model is a strict conditional refinement of first-order MAO theory and preserves observed pairwise association while allowing higher-order membership structure to vary.

Using membership atoms, we derived exact formulas for the partition function, the complete distributions of exact-$t$ and at-least-$t$ occupancy counts, and arbitrary factorial moments. The second-order MAO norms satisfy

$$\mathbb{E}[(X_{=t})_r] = \|t^r\|_T,$$

and

$$\mathbb{E}[(X_{\geq t})_r] = \|[t, T]^r\|_T.$$

Thus, all ordinary moments, central moments, and cumulants of both occupancy families are available from a common exact framework.

The validation results confirm exact agreement among atom-weighted PMFs, direct atomic factorial-moment calculations, and brute-force labelled-set enumeration where feasible. The threshold statistic $X_{\geq t}$ substantially broadens the practical utility of the theory because it directly represents repeated, persistent, or widely shared occurrence.

The natural next step is higher-order MAO theory, in which intersections of order three and above are also fixed. The atomic formulation persists, but scalable computation and asymptotic analysis will require new methods.

**Acknowledgements**

An OpenAI language model (ChatGPT) was used for language editing and limited assistance with data presentation and visualization. The author independently developed the underlying mathematical theory and its core computational implementation, independently validated the results, and takes full responsibility for the manuscript.

**Appendix A. First-order MAO norms as mixtures of second-order MAO norms**

Let

$$[T] = \{1, \dots, T\},$$

and let

$$m = (m_1, \dots, m_T),\ 0 \le m_i \le N.$$

The first-order MAO configuration space is

$$\Omega_1(N; m) = \{(A_1, \dots, A_T): A_i \subseteq [N]\ \ |A_i| = m_i\ \ i \in [T]\}. \qquad \text{(A.1)}$$

Under the first-order MAO model, $\Omega_1(N; m)$ is equipped with the uniform distribution. Equivalently, the sets $A_1, \dots, A_T$ are independently and uniformly sampled subject to their prescribed cardinalities.

For a configuration $(A_1, \dots, A_T)$, define its pairwise-overlap matrix by

$$C = (c_{ij})_{1 \le i,j \le T},\ c_{ii} = m_i,\ c_{ij} = |\ A_i \cap A_j\ |\ \ (i \ne j). \qquad \text{(A.2)}$$

Let $\mathcal{C}(m)$ denote the set of feasible pairwise-overlap matrices compatible with $N$ and $m$. For each $C \in \mathcal{C}(m)$, let

$$\Omega_2(N; (m, C)) = \{(A_1, \dots, A_T) \in \Omega_1(N; m):|\ A_i \cap A_j\ | = c_{ij} \text{ for all } i, j \in [T]\} \qquad \text{(A.3)}$$

The first-order configuration space is partitioned into disjoint second-order fibers:

$$\Omega_1(N; m) = \bigsqcup_{C \in \mathcal{C}(m)} \Omega_2\,(N; (m, C)). \qquad \text{(A.4)}$$

Consequently, if

$$P_{N,m}(C) = \frac{|\ \Omega_2(N; (m, C))\ |}{|\ \Omega_1(N; m)\ |}, \qquad \text{(A.5)}$$

then, for every statistic $Y$ defined on $\Omega_1(N; m)$,

$$\mathbb{E}_{N,m}[Y] = \sum_{C \in \mathcal{C}(m)} P_{N,m}\,(C)\mathbb{E}_{N,(m,C)}[Y]. \qquad \text{(A.6)}$$

The purpose of this appendix is to show that, when $Y$ is a factorial occupancy statistic, the mixture identity (A.6) recovers the original first-order MAO norm formula.

**A.1 First-order joint membership probabilities**

Fix $r \ge 1$, and let

$$u_1, \dots, u_r \in [N]$$

be distinct labeled individuals. For each $j \in [r]$, let

$$B_j \subseteq [T]$$

denote the membership pattern of individual $u_j$, so that

$$u_j \in A_i \iff i \in B_j. \qquad \text{(A.7)}$$

For every set index $i \in [T]$, define

$$k_i(B_1, \dots, B_r) = \sum_{j=1}^{r} 1\{i \in B_j\}. \qquad \text{(A.8)}$$

Thus, $k_i(B_1, \dots, B_r)$ is the number of the $r$ labeled individuals that are required to belong to $A_i$.

Define the first-order MAO allocation function

$$g_m(B_1, \dots, B_r) = \prod_{i=1}^{T} (m_i)_{k_i(B_1,\dots,B_r)} (N - m_i)_{r - k_i(B_1,\dots,B_r)}. \qquad \text{(A.9)}$$

For one set $A_i$, sampled uniformly from all $m_i$-subsets of $[N]$, the probability of the required membership pattern is

$$\mathbb{P}_{N,m}\big(u_j \in A_i \Leftrightarrow i \in B_j \ 1 \le j \le r\big) = \frac{(m_i)_{k_i(B_1,\dots,B_r)} (N - m_i)_{r - k_i(B_1,\dots,B_r)}}{(N)_r}. \qquad \text{(A.10)}$$

Since $A_1, \dots, A_T$ are independent under the first-order MAO model,

$$\begin{aligned} &\mathbb{P}_{N,m}\big(\{i \in [T]: u_j \in A_i\} = B_j \ 1 \le j \le r\big) \\ &= \prod_{i=1}^{T} \frac{(m_i)_{k_i(B_1,\dots,B_r)} (N - m_i)_{r - k_i(B_1,\dots,B_r)}}{(N)_r} \\ &= \frac{g_m(B_1, \dots, B_r)}{(N)_r^T}. \end{aligned} \qquad \text{(A.11)}$$

Equation (A.11) is the probabilistic origin of the first-order MAO allocation function $g_m$.

### A.2 First-order exact-occupancy norms

Recall that

$$X_{=t} = \sum_{u=1}^{N} 1\{|\{i \in [T]: u \in A_i\}| = t\}. \qquad \text{(A.12)}$$

The falling factorial $(X_{=t})_r$ counts ordered $r$-tuples of distinct individuals having occupancy exactly $t$. By exchangeability of the $N$ individual labels,

$$\mathbb{E}_{N,m}[(X_{=t})_r] = (N)_r \mathbb{P}_{N,m}(|B_j| = t\ 1 \le j \le r). \qquad \text{(A.13)}$$

Using (A.11),

$$\begin{aligned}\mathbb{E}_{N,m}[(X_{=t})_r] &= (N)_r \sum_{\substack{B_1,\dots,B_r \subseteq [T] \\ |B_1| = \cdots = |B_r| = t}} \frac{g_m(B_1,\dots,B_r)}{(N)_r^T} \\ &= \frac{\sum_{\substack{B_1,\dots,B_r \subseteq [T] \\ |B_1| = \cdots = |B_r| = t}} g_m(B_1,\dots,B_r)}{(N)_r^{T-1}}.\end{aligned} \qquad \text{(A.14)}$$

Because the first-order exact-$t$ MAO norm is defined by

$$\| t^r \|_{T,1} = \mathbb{E}_{N,m}[(X_{=t})_r], \qquad \text{(A.15)}$$

we obtain the original first-order MAO norm formula:

$$\boxed{\| t^r \|_{T,1} = \frac{\sum_{\substack{B_1,\dots,B_r \subseteq [T] \\ |B_1| = \cdots = |B_r| = t}} g_m(B_1,\dots,B_r)}{(N)_r^{T-1}}.} \qquad \text{(A.16)}$$

This is precisely the earlier first-order MAO norm expression, now derived from the factorial-moment interpretation.

### A.3 First-order threshold-occupancy norms

Similarly,

$$X_{\ge t} = \sum_{u=1}^{N} 1\{|\{i \in [T]: u \in A_i\}| \ge t\}. \qquad \text{(A.17)}$$

Again, $(X_{\ge t})_r$ counts ordered $r$-tuples of distinct individuals whose occupancy is at least $t$. Therefore,

$$\mathbb{E}_{N,m}[(X_{\ge t})_r] = (N)_r \mathbb{P}_{N,m}(|B_j| \ge t\ 1 \le j \le r). \qquad \text{(A.18)}$$

Applying (A.11) yields

$$\mathbb{E}_{N,m}[(X_{\geq t})_r] = \frac{\sum_{\substack{B_1,\dots,B_r \subseteq [T] \\ |B_j| \geq t,\ 1 \leq j \leq r}} g_m(B_1,\dots,B_r)}{(N)_r^{T-1}}. \tag{A.19}$$

Hence,

$$\boxed{\|[t,T]^r\|_{T,1} = \frac{\sum_{\substack{B_1,\dots,B_r \subseteq [T] \\ |B_j| \geq t,\ 1 \leq j \leq r}} g_m(B_1,\dots,B_r)}{(N)_r^{T-1}}.} \tag{A.20}$$

Equation (A.20) recovers the original first-order threshold MAO norm formula.

### A.4 Mixture representation through second-order MAO norms

For each feasible pairwise-overlap matrix $C \in \mathcal{C}(m)$, the conditional distribution of

$$(A_1, \dots, A_T)$$

given the event that its pairwise-overlap matrix equals $C$ is uniform over

$$\Omega_2(N; (m, C)). \tag{A.21}$$

Thus, for every statistic $Y$,

$$\mathbb{E}_{N,m}[Y \mid C] = \mathbb{E}_{N,(m,C)}[Y]. \tag{A.22}$$

Taking

$$Y = (X_{=t})_r$$

in (A.6) gives

$$\begin{aligned} \| t^r \|_{T,1} &= \mathbb{E}_{N,m}[(X_{=t})_r] \\ &= \sum_{C \in \mathcal{C}(m)} P_{N,m}(C)\, \mathbb{E}_{N,(m,C)}[(X_{=t})_r] \\ &= \sum_{C \in \mathcal{C}(m)} P_{N,m}(C)\, \| t^r \|_{T,2;(m,C)}. \end{aligned} \tag{A.23}$$

Likewise, taking

$$Y = (X_{\geq t})_r$$

gives

$$\begin{aligned}\|[t,T]^r\|_{T,1} &= \mathbb{E}_{N,m}[(X_{\geq t})_r] \\ &= \sum_{C\in\mathcal{C}(m)} P_{N,m}(C)\mathbb{E}_{N,(m,C)}[(X_{\geq t})_r] \\ &= \sum_{C\in\mathcal{C}(m)} P_{N,m}(C)\|[t,T]^r\|_{T,2;(m,C)}.\end{aligned} \tag{A.24}$$

Combining (A.16) with (A.23) gives

$$\boxed{\frac{\sum_{\substack{B_1,\dots,B_r\subseteq[T]\\|B_1|=\cdots=|B_r|=t}} g_m(B_1,\dots,B_r)}{(N)_r^{T-1}} = \sum_{C\in\mathcal{C}(m)} P_{N,m}(C)\ \|\ t^r\ \|_{T,2;(m,C)}.} \tag{A.25}$$

Similarly, combining (A.20) with (A.24) gives

$$\boxed{\frac{\sum_{\substack{B_1,\dots,B_r\subseteq[T]\\|B_j|\geq t,\ 1\leq j\leq r}} g_m(B_1,\dots,B_r)}{(N)_r^{T-1}} = \sum_{C\in\mathcal{C}(m)} P_{N,m}(C)\|[t,T]^r\|_{T,2;(m,C)}.} \tag{A.26}$$

Equations (A.25) and (A.26) are the desired reduction formulas. They show that the original first-order MAO norm is exactly the mixture, over the pairwise-overlap structures induced by first-order conditioning, of the corresponding second-order MAO norms.

In particular, the second-order MAO norm construction extends rather than replaces the first-order MAO norm construction:

$$\boxed{\text{first-order MAO norm} = \text{mixture of second-order MAO norms.}}$$

The explicit first-order allocation formula involving $g_m$ is recovered after the second-order fibers are aggregated according to their first-order probabilities $P_{N,m}(C)$.

**A.5 Summary: recovery of the first-order MAO norm**

The preceding argument establishes two equivalent representations of the same factorial moment. For every fixed $r \geq 1$,

$$\boxed{\|\ t^r\ \|_{T,1} = \underbrace{\frac{\sum_{\substack{B_1,\dots,B_r\subseteq[T]\\|B_1|=\cdots=|B_r|=t}} g_m(B_1,\dots,B_r)}{(N)_r^{T-1}}}_{\text{direct first-order MAO calculation}} = \underbrace{\sum_{C\in\mathcal{C}(m)} P_{N,m}(C)\ \|\ t^r\ \|_{T,2;(m,C)}}_{\text{mixture of second-order MAO norms}}.} \tag{A.27}$$

The left-hand expression is obtained directly from the first-order MAO model. Independent uniform sampling of $A_1,\dots,A_T$, conditional only on the set sizes $m_1,\dots,m_T$, gives the joint membership probability

$$\mathbb{P}_{N,m}(B_1, \dots, B_r) = \frac{g_m(B_1, \dots, B_r)}{(N)_r^T}.$$

Multiplying by the number $(N)_r$ of ordered $r$-tuples of distinct individuals and summing over all membership patterns of size $t$ yields the original first-order MAO norm formula.

The right-hand expression is obtained by partitioning the first-order configuration space according to its induced pairwise-overlap matrix:

$$\Omega_1(N; m) = \bigsqcup_{C \in \mathcal{C}(m)} \Omega_2(N; (m, C)).$$

Within every fiber $\Omega_2(N; (m, C))$, the conditional first-order distribution is exactly the uniform second-order MAO distribution. Therefore, the law of total expectation gives the weighted average of the corresponding second-order factorial moments.

Thus, aggregating the second-order MAO norms over all feasible pairwise-overlap structures reproduces the first-order MAO norm exactly. The result is not merely a decomposition of the first-order configuration space: it shows that the second-order MAO norm construction is consistent with, and recovers, the original first-order allocation formula based on $g_m$.

For threshold occupancy, the same argument gives

$$\boxed{\|[t, T]^r\|_{T,1} = \frac{\sum_{\substack{B_1, \dots, B_r \subseteq [T] \\ |B_j| \geq t,\ 1 \leq j \leq r}} g_m(B_1, \dots, B_r)}{(N)_r^{T-1}} = \sum_{C \in \mathcal{C}(m)} P_{N,m}(C) \|[t, T]^r\|_{T,2;(m,C)}.} \tag{A.28}$$

**Appendix B. Poisson convergence via factorial moments**

Let $\{M^{(N)}\}_{N\geq 1}$ be a feasible sequence of second-order MAO constraint matrices. For each $N$, let $\mathbb{P}_{N,M^{(N)}}$ denote the uniform probability measure on the corresponding feasible second-order MAO configurations. Let $X_N$ denote either $X_{=t,N}$ or $X_{\geq t,N}$, where $t$ is fixed.

**Theorem B.1.** Suppose that, for every fixed integer $r \geq 1$,

$$\mathbb{E}_{N,M^{(N)}}\,[(X_N)_r] \longrightarrow \lambda^r \qquad \text{(B.1)}$$

for some $\lambda \in [0,\infty)$. Then

$$X_N \Rightarrow \text{Poisson}(\lambda). \qquad \text{(B.2)}$$

For $X_N = X_{=t,N}$, condition (B.1) is equivalently

$$\|t^r\|^{(2)}_{T,M^{(N)}} \longrightarrow \lambda^r \text{ for every fixed } r \geq 1. \qquad \text{(B.3)}$$

For $X_N = X_{\geq t,N}$, condition (B.1) is equivalently

$$\|[t,T]^r\|^{(2)}_{T,M^{(N)}} \longrightarrow \lambda^r \text{ for every fixed } r \geq 1 \qquad \text{(B.4)}$$

**Proof.** Let $Z \sim \text{Poisson}(\lambda)$. Its factorial moments satisfy

$$\mathbb{E}[(Z)_r] = \lambda^r \text{ for every } r \geq 0, \qquad \text{(B.5)}$$

where $(Z)_0 = 1$. Since $(X_N)_0 = 1$ almost surely, (B.1) and (B.5) give

$$\mathbb{E}_{N,M^{(N)}}\,[(X_N)_r] \longrightarrow \mathbb{E}[(Z)_r] \text{ for every fixed } r \geq 0. \qquad \text{(B.6)}$$

For every integer $v \geq 1$, the power $x^v$ has the factorial expansion

$$x^v = \sum_{i=1}^{v} s_{v,i}\,(x)_i, \qquad \text{(B.7)}$$

where $s_{v,i}$ denotes a Stirling number of the second kind. Therefore,

$$\mathbb{E}_{N,M^{(N)}}[X_N^v] = \sum_{i=1}^{v} s_{v,i}\,\mathbb{E}_{N,M^{(N)}}\,[(X_N)_i]. \qquad \text{(B.8)}$$

Because the sum in (B.8) is finite, (B.6) yields

$$
\begin{aligned}
\mathbb{E}_{N,M^{(N)}}[X_N^v] &\longrightarrow \sum_{i=1}^{v} s_{v,i}\, \mathbb{E}[(Z)_i] \\
&= \sum_{i=1}^{v} s_{v,i}\, \lambda^i \\
&= \mathbb{E}[Z^v].
\end{aligned}
\tag{B.9}
$$

Thus, for every fixed $v \geq 1$, the $v$-th ordinary moment of $X_N$ converges to the corresponding ordinary moment of $Z$. Since the Poisson distribution is moment-determinate, the method of moments implies

$$X_N \Rightarrow Z \sim \text{Poisson}(\lambda).$$

This proves (B.2).

The proof is identical to the first-order factorial-moment argument once (B.1) has been established. The distinction lies in verifying (B.1): under second-order constraints, the prescribed pairwise intersections restrict the full symmetric first-order configuration class, so the required factorial-moment asymptotics cannot in general be inferred from the first-order product structure alone.

## Appendix C. Conditional Central Limit Lemma

This appendix supplies the probabilistic step used in Section 7.2. We begin with independent membership vectors, for which an occupancy count is an ordinary sum of iid indicators. The uniform MAO law is obtained by conditioning that **independent and identically distributed (iid)** sample on having exactly the prescribed first-order and second-order totals. The lemma shows that, after fixing the set sizes and pairwise intersections, the occupancy count is still asymptotically normal if it is not completely determined by those fixed quantities.

**Informal interpretation.** We work only with regular second-order constraint matrices: the prescribed set sizes and pairwise intersections must be attainable, must lie away from the boundary of the feasible region, and must leave genuine freedom in the higher-order membership patterns. Under these conditions, fixing the first- and second-order totals removes some variation in an occupancy count, but does not determine it completely. The remaining fluctuations are of order $\sqrt{n}$ and are asymptotically normal.

### C.1 The conditional-sum setting

We describe one element at a time. Each element matters in two ways. First, it may contribute to the occupancy count that we want to study. Second, it contributes to the set sizes and pairwise intersections that are held fixed in the MAO model.

For element $r$, let $H_r$ record its contribution to the occupancy count, and let $U_r$ record its contribution to the fixed constraints. We keep these two quantities together because they are produced by the same element and may be correlated.

For example, for the exact-$t$ occupancy count, $H_r = 1$ if element $r$ belongs to exactly $t$ sets, and $H_r = 0$ otherwise. The vector $U_r$ records which set sizes and pairwise intersections receive a contribution of one from that same element.

Adding over all $n$ elements gives

$$W_n := \sum_{r=1}^{n} H_r \, , \; S_n := \sum_{r=1}^{n} U_r \, . \qquad \text{(C.1)}$$

Thus $W_n$ is the occupancy count of interest (target contribution), while $S_n$ is the vector of set sizes and pairwise intersections that will be fixed exactly (constraint contribution). Under the iid reference model, $S_n$ is the random vector of total set sizes and pairwise intersections. We choose the reference distribution so that its mean matches the prescribed MAO constraints. Thus, if $\mu = \mathbb{E}[U_1]$, the prescribed total is $n\mu$. We then condition on the event $S_n = n\mu$, meaning that every set size and every pairwise intersection has exactly its prescribed value.

There are three technical requirements. First, the possible values of $U_r$ must genuinely vary in every retained constraint direction; otherwise redundant directions should be removed. Second, the exact value $n\mu$ must be attainable along the sequence considered. Third, the target must not be completely determined by the constraints. The last requirement is quantified by the residual variance $v$ below.

Let $L$ be the lattice generated by differences of support points of $U_1$. We assume that, after redundant coordinates have been removed, $L$ has full rank and $U_1$ is aperiodic on this lattice. This is the standard lattice hypothesis which allows a local limit theorem at the exact conditioning point.

### C.2 The lemma

**Lemma C.1 (Conditional CLT at an exact lattice constraint).** Let $(H_1, U_1)_{r\geq 1}$ be iid, with $H_r \in \mathbb{R}$ and $U_r \in \mathbb{Z}^d$. Assume that $(H_1, U_1)$ has finite support. Set

$$\mu := \mathbb{E}[U_1],\ m := \mathbb{E}[H_1],\ \Sigma_U := \mathrm{Var}(U_1). \qquad \text{(C.2)}$$

Assume that $\Sigma_U$ is positive definite and that $U_1$ is strongly aperiodic on its minimal supporting lattice. Let $\mathcal{N} \subseteq \mathbb{N}$ be an infinite set for which

$$\mathbb{P}(S_n = n\mu) > 0\ (n \in \mathcal{N}). \qquad \text{(C.3)}$$

The covariance between the target contribution and the constraint contribution is

$$c := \mathrm{Cov}(U_1, H_1). \qquad \text{(C.4)}$$

The part of the variance of $H_1$ that remains after linearly accounting for $U_1$ is

$$v := \mathrm{Var}(H_1) - c^T \Sigma_U^{-1} c. \qquad \text{(C.5)}$$

If $v > 0$, then, as $n \to \infty$ through $n \in \mathcal{N}$,

$$\mathcal{L}\left(\frac{W_n - \mathbb{E}[W_n \mid S_n = n\mu]}{\sqrt{nv}}\middle| S_n = n\mu\right) \Rightarrow N(0,1), \qquad \text{(C.6)}$$

and

$$\mathrm{Var}(W_n \mid S_n = n\mu) = nv + o(n). \qquad \text{(C.7)}$$

### C.3 Why the coefficient is $v$

The vector $\beta := \Sigma_U^{-1} c$ is the coefficient of the linear least-squares prediction of $H_1$ from $U_1$. Define the unpredictable remainder by

$$R_r := H_r - m - \beta^T (U_r - \mu). \qquad \text{(C.8)}$$

By the normal equations for linear projection,

$$\mathbb{E}[R_r] = 0,\ \mathrm{Cov}(R_r, U_r) = 0,\ \mathrm{Var}(R_r) = v. \tag{C.9}$$

Summing (C.8) gives the exact identity

$$W_n - nm = \beta^T(S_n - n\mu) + \sum_{r=1}^{n} R_r\,. \tag{C.10}$$

On the conditioning event $S_n = n\mu$, the first term vanishes. Hence

$$W_n - nm = \sum_{r=1}^{n} R_r \text{ when } S_n = n\mu. \tag{C.11}$$

This identity is the whole reason for the Schur-complement formula: after the totals have been fixed, only the residuals can fluctuate.

**C.4 Proof of Lemma C.1**

Let $D$ be a fundamental domain of the dual lattice. Fourier inversion expresses the conditional characteristic function of the residual sum as

$$\mathbb{E}\left[\exp\left\{\frac{is}{\sqrt{n}}\sum_{r=1}^{n} R_r\right\}\middle| S_n = n\mu\right] = \frac{\int_D \left[\mathbb{E}e^{isR_1/\sqrt{n}+i\theta^T(U_1-\mu)}\right]^n d\theta}{\int_D \left[\mathbb{E}e^{i\theta^T(U_1-\mu)}\right]^n d\theta}. \tag{C.12}$$

Near $\theta = 0$, put $\theta = u/\sqrt{n}$. Finite support gives a uniform Taylor expansion, and the orthogonality in (C.9) removes the mixed term:

$$\mathbb{E}e^{isR_1/\sqrt{n}+iu^T(U_1-\mu)/\sqrt{n}} = 1 - \frac{vs^2 + u^T\Sigma_U u}{2n} + O\left(n^{-3/2}\right). \tag{C.13}$$

Consequently,

$$\left[\mathbb{E}e^{isR_1/\sqrt{n}+iu^T(U_1-\mu)/\sqrt{n}}\right]^n \longrightarrow \exp\left\{-\frac{1}{2}vs^2 - \frac{1}{2}u^T\Sigma_U u\right\}. \tag{C.14}$$

Aperiodicity makes the portions of the integrals in (C.12) away from $\theta = 0$ exponentially small. On a fixed neighbourhood of zero, (C.13), the substitution $\theta = u/\sqrt{n}$, and dominated convergence give the usual multivariate lattice local-limit calculation. The Gaussian integral involving $u$ occurs in both numerator and denominator and cancels. Thus

$$\mathbb{E}\left[\exp\left\{\frac{is}{\sqrt{n}}\sum_{r=1}^{n} R_r\right\}\middle| S_n = n\mu\right] \longrightarrow e^{-vs^2/2}. \tag{C.15}$$

Levy's theorem yields

$$\mathcal{L}\left(\frac{1}{\sqrt{n}}\sum_{r=1}^{n} R_r \middle| S_n = n\mu\right) \Rightarrow N(0, v). \qquad \text{(C.16)}$$

The same local-limit expansion, differentiated at $s = 0$, gives

$$\frac{1}{\sqrt{n}}\mathbb{E}\left[\sum_{r=1}^{n} R_r \middle| S_n = n\mu\right] \to 0,\ \frac{1}{n}\mathbb{E}\left[\left(\sum_{r=1}^{n} R_r\right)^2 \middle| S_n = n\mu\right] \to v. \qquad \text{(C.17)}$$

Using (C.11), these two limits respectively replace $nm$ by the conditional mean in (C.16), and yield (C.7). This proves the lemma.

**C.5 Translation to the second-order MAO model**

We now identify the abstract objects with the model of Section 7.2. A membership vector $Y_u \in \{0, 1\}^T$ records which of the $T$ sets contain element $u$. One element contributes 1 to a set size $M_{ii}$ when $Y_{u,i} = 1$, and contributes 1 to an intersection $M_{ij}$ when $Y_{u,i}Y_{u,j} = 1$. Thus define

$$\Psi(y) := \big(y_i, , y_i y_j : i \in [T]\ 1 \le i < j \le T\big). \qquad \text{(C.18)}$$

After retaining a maximal linearly independent collection of these coordinates, the abstract constraint variable is

$$U_r = \Psi(Y_r^\star). \qquad \text{(C.19)}$$

Its sum is precisely the prescribed vector of set sizes and pairwise intersections. The target contribution is simply the indicator that one element has the desired occupancy:

$$H_r = 1\left\{\sum_{i=1}^{T} Y_{r,i}^\star = t\right\} \text{ for } X_{=t}, \qquad \text{(C.20)}$$

or

$$H_r = 1\left\{\sum_{i=1}^{T} Y_{r,i}^\star \ge t\right\} \text{ for } X_{\ge t}. \qquad \text{(C.21)}$$

For $N_q = qN_0$ and $M^{(N_q)} = qM^{(0)}$, the required total constraint is attainable for every $q$.

The maximum-entropy iid law, conditioned on that exact total, is uniform on the corresponding MAO fibre, because its probability depends only on that total.

Accordingly, Section 7.2 needs only two substantive assumptions:

- **Interior constraints.** The normalized sizes and intersections lie in the relative interior of their feasible moment region. In concrete terms, no membership pattern is asymptotically forced to disappear. This supplies the nondegenerate iid reference law required by the local-limit argument.

- **Unforced statistic.** After the sizes and pairwise intersections have been fixed, the chosen occupancy count still has random variation. In formulae this is exactly $v > 0$, with $v$ defined in (C.5).

All other conditions in Lemma C.1 are bookkeeping consequences of this setup: finite support follows from fixed $T$; exact attainability follows from feasible proportional scaling; and redundant constraint coordinates are removed once before the lemma is applied.

**Supplementary Material 1**

**Supplementary Note: Membership Atoms, Labelled Configurations, and the Partition Function**

This note gives an interpretation of the objects appearing in equations (1)--(10). The main purpose is to distinguish: (i) the membership pattern of one labelled element, (ii) the vector of membership-atom counts, and (iii) a complete labelled set system.

Let $[T] = \{1, \dots, T\}$, let $[N] = \{1, .., N\}$, and let $A_1, \dots, A_T \subseteq [N]$ be labelled sets.

**Membership patterns and atoms**

Each labelled element $u \in [N]$ has a binary membership vector

$$Y_u = (1\{u \in A_1\}, \dots, 1\{u \in A_T\}) \in \{0,1\}^T.$$

Equivalently, $u$ has an exact membership set

$$B_u = \{i \in [T]: u \in A_i\} \subseteq [T].$$

For every $B \subseteq [T]$, define the membership atom

$$C_B = \left(\bigcap_{i \in B} A_i\right) \cap \left(\bigcap_{j \notin B} A_j^c\right), \; a_B = | C_B |.$$

Thus, $a_B$ is the number of elements belonging to exactly the sets indexed by $B$, and to no other sets. The $2^T$ atoms are mutually disjoint and partition the population:

$$[N] = \bigsqcup_{B \subseteq [T]} C_B.$$

Consequently,

$$\sum_{B \subseteq [T]} a_B = N. \qquad (4)$$

Equation (4) merely states that every element has exactly one complete membership pattern. An element belongs to one atom, but may belong to several original sets. For example, an element in $C_{\{1,3\}}$ belongs to $A_1$ and $A_3$, but not to any other set.

For $T = 3$, the atom-count vector is

$$a = \left(a_\emptyset, a_{\{1\}}, a_{\{2\}}, a_{\{3\}}, a_{\{1,2\}}, a_{\{1,3\}}, a_{\{2,3\}}, a_{\{1,2,3\}}\right).$$

These eight entries are precisely the eight cell counts of the three-set Venn diagram. In particular,

$$a_{\{1,2,3\}} = | A_1 \cap A_2 \cap A_3 |$$

is the count in the central three-way overlap, whereas

$$a_{\{1,2\}}$$

counts elements in $A_1 \cap A_2$ that are not in $A_3$.

**Translation of marginal and pairwise constraints**

The set size $| A_i |$ is obtained by adding all atoms whose membership pattern contains $i$:

$$\sum_{B:\, i \in B} a_B = M_{ii},\ i \in [T]. \qquad (5)$$

For example, when $T = 3$,

$$| A_1 | = a_{\{1\}} + a_{\{1,2\}} + a_{\{1,3\}} + a_{\{1,2,3\}}.$$

The pairwise intersection $| A_i \cap A_j |$ is obtained by adding all atoms whose membership pattern contains both $i$ and $j$:

$$\sum_{B:\, \{i,j\} \subseteq B} a_B = M_{ij},\ 1 \leq i < j \leq T. \qquad (6)$$

For $T = 3$,

$$| A_1 \cap A_2 | = a_{\{1,2\}} + a_{\{1,2,3\}}.$$

The second term is essential: every element in the three-way intersection also belongs to $A_1 \cap A_2$. More generally, for $T = 4$,

$$| A_1 \cap A_2 | = a_{\{1,2\}} + a_{\{1,2,3\}} + a_{\{1,2,4\}} + a_{\{1,2,3,4\}}.$$

Thus, equations (5) and (6) translate the observed marginal set sizes and pairwise intersections into linear constraints on the Venn-diagram cell counts.

**Three distinct levels of description**

The following objects should be kept distinct:

element-level membership vectors $\longrightarrow$ atom-count vectors $\longrightarrow$ labelled set systems.

| Object | Form | Meaning |
|---|---|---|
| $Y_u$ | $\{0,1\}^T$ | Membership pattern of one labelled element $u$ |
| $a = (a_B)_{B\subseteq[T]}$ | $\mathbb{Z}_{\geq 0}^{2^T}$ | Counts of elements having each membership pattern |
| $(A_1, \dots, A_T)$ | $T$ labelled subsets of $[N]$ | One complete labelled configuration |
| $\mathcal{A}_2(N, M)$ | a set of atom vectors | All feasible Venn-diagram count configurations |
| $\Omega_2(N, M)$ | a set of labelled set systems | All feasible labelled configurations |

In particular,

$$\mathcal{A}_2(N, M) = \left\{(a_B)_{B\subseteq[T]} \in \mathbb{Z}_{\geq 0}^{2^T}\text{:(4), (5), and (6) hold}\right\}.$$

Hence, $\mathcal{A}_2(N, M)$ is not a collection of original set systems and is not a set of binary vectors. It is the set of all feasible **atom-count vectors**, or equivalently, all feasible numerical Venn diagrams.

The constraints generally do not determine a unique $a$. For example, even when all set sizes and all pairwise intersections are fixed, the count

$$a_{\{1,2,3\}} = | A_1 \cap A_2 \cap A_3 |$$

need not be fixed. This remaining freedom is precisely the higher-order membership structure randomized by the second-order MAO model.

**From atom counts to labelled configurations**

A feasible atom vector $a$ specifies how many elements occupy each atom, but does not specify which labelled elements occupy those atoms. For example,

$$a_{\{1,2,3\}} = r$$

states that exactly $r$ elements belong to all three sets, but does not identify those elements.

Given $a$, the $N$ labelled elements must be allocated to the $2^T$ labelled atoms, with atom $B$ receiving exactly $a_B$ elements. The number of such allocations is the multinomial coefficient

$$w(a) = \frac{N!}{\prod_{B\subseteq[T]} a_B\,!}. \qquad (7)$$

This is the standard count of allocations of $N$ labelled objects to labelled boxes with prescribed occupancies. The factor $a_B!$ removes permutations of elements within atom $B$, since their order within that atom is irrelevant.

Once every element has been assigned to an atom, the sets are uniquely determined:

$$A_i = \{u \in [N]: u \text{ is assigned to an atom } B \text{ with } i \in B\}.$$

Conversely, every labelled set system determines a unique allocation of its elements to the membership atoms. Therefore, labelled atom allocations and labelled set systems are in one-to-one correspondence.

**Configuration space and partition function**

Define

$$\Omega_2(N, M) = \{(A_1, \dots, A_T): \text{the prescribed set sizes and pairwise intersections hold}\}.$$

The space $\Omega_2(N, M)$ is the set of all concrete labelled set systems satisfying the second-order constraints. It is naturally partitioned according to the atom vector that each configuration induces:

$$\Omega_2(N, M) = \bigsqcup_{a \in \mathcal{A}_2(N, M)} \Omega\,(a),$$

where

$$\Omega(a) = \{(A_1, \dots, A_T) \in \Omega_2(N, M): \text{the induced atom vector is } a\}.$$

For each feasible atom vector,

$$|\,\Omega(a)\,| = w(a).$$

The partition function is therefore

$$Z_{N,M}^{(2)} = \sum_{a \in \mathcal{A}_2(N, M)} \frac{N!}{\prod_{B \subseteq [T]} a_B\,!}. \qquad (8)$$

In the present setting, the term "partition function" has a direct combinatorial meaning: it is the sum of the sizes of all atom-vector classes. Hence,

$$Z_{N,M}^{(2)} = \sum_{a \in \mathcal{A}_2(N, M)} |\,\Omega(a)\,| = |\,\Omega_2(N, M)\,|. \qquad (9)$$

Thus, $Z_{N,M}^{(2)}$ is both:

- the total number of labelled set systems satisfying the constraints; and

- the normalizing constant required to define probabilities on that constrained configuration space.

**The induced atom-vector distribution**

Under the second-order MAO null model, every labelled configuration in $\Omega_2(N, M)$ is equally likely. Therefore,

$$\mathbb{P}^{(2)}_{N,M}(\omega) = \frac{1}{Z^{(2)}_{N,M}}, \ \omega \in \Omega_2(N, M).$$

After the element labels are ignored and only the atom vector is retained, the induced probability of a feasible vector $a$ is

$$\mathbb{P}^{(2)}_{N,M}(a) = \frac{w(a)}{Z^{(2)}_{N,M}} = \frac{N!/\prod_{B\subseteq[T]} a_B\,!}{Z^{(2)}_{N,M}}, \ a \in \mathcal{A}_2(N, M). \tag{10}$$

This distribution is generally not uniform over $\mathcal{A}_2(N, M)$. Different atom vectors may correspond to different numbers of labelled set systems. A feasible Venn-diagram count configuration with a larger multinomial weight represents more labelled allocations and therefore has larger probability under a null model that is uniform over labelled configurations.

**Recovery of MAO occupancy statistics**

Finally, the exact and threshold occupancy counts are functions of the atom vector:

$$X_{=t} = \sum_{\substack{B\subseteq[T] \\ |B|=t}} a_B \,, \ t = 0, \ldots, T,$$

and

$$X_{\geq t} = \sum_{\substack{B\subseteq[T] \\ |B|\geq t}} a_B = \sum_{s=t}^{T} X_{=s}.$$

For $T = 3$, these identities become

$$X_{=0} = a_{\emptyset},$$

$$X_{=1} = a_{\{1\}} + a_{\{2\}} + a_{\{3\}},$$

$$X_{=2} = a_{\{1,2\}} + a_{\{1,3\}} + a_{\{2,3\}},$$

and

$$X_{=3} = a_{\{1,2,3\}}.$$

Therefore, the atom-vector formulation is a detailed representation of membership structure, while $X_{=t}$ and $X_{\geq t}$ are its occupancy summaries. Equations (1)--(3) express the deterministic first- and second-order constraints on those summaries; equations (4)--(10) provide the exact constrained distribution from which their remaining higher-order variation is derived.

**A concrete three-set example: configurations, atom vectors, and the partition function**

Let the labelled ground set be

$$U = \{e_1, e_2, e_3, e_4, e_5\}, \quad |U| = N = 5.$$

Consider three labelled sets

$$A_1, A_2, A_3 \subseteq U.$$

Suppose that their marginal sizes and pairwise intersections are prescribed by

$$|A_1| = |A_2| = |A_3| = 2,$$

and

$$|A_1 \cap A_2| = |A_1 \cap A_3| = |A_2 \cap A_3| = 1.$$

Equivalently, the prescribed second-order overlap matrix is

$$M = \begin{pmatrix} 2 & 1 & 1 \\ 1 & 2 & 1 \\ 1 & 1 & 2 \end{pmatrix}.$$

The entries of $M$ have a direct meaning:

$$M_{ii} = |A_i| = 2, \; M_{ij} = |A_i \cap A_j| = 1 \; (i \neq j).$$

The important point is that these constraints do **not** determine the triple intersection

$$|A_1 \cap A_2 \cap A_3|.$$

There are two possible higher-order overlap patterns.

**1. The configuration space $\Omega_2(U, M)$**

For this particular ground set $U$, define

$$\Omega_2(U,M)=\left\{(A_1,A_2,A_3):\begin{array}{c} A_i\subseteq U,\\ |A_i|=2\ (i=1,2,3),\\ |A_i\cap A_j|=1\ (1\le i<j\le 3)\end{array}\right\}.$$

An element of $\Omega_2(U,M)$ is a concrete labelled triple of subsets of $U$.

For example, consider

$$\omega^{(0)}=(A_1,A_2,A_3)=(\{e_1,e_2\},\{e_1,e_3\},\{e_2,e_3\}).$$

Here,

$$A_1\cap A_2=\{e_1\},\ A_1\cap A_3=\{e_2\},\ A_2\cap A_3=\{e_3\}.$$

The three pairwise intersections are realized by three distinct elements. Moreover,

$$A_1\cap A_2\cap A_3=\emptyset.$$

The elements $e_4,e_5$ belong to none of the three sets.

Now consider the different configuration

$$\omega^{(1)}=(A_1,A_2,A_3)=(\{e_1,e_2\},\{e_1,e_3\},\{e_1,e_4\}).$$

Again,

$$|A_1|=|A_2|=|A_3|=2,$$

and all pairwise intersections have size one:

$$A_1\cap A_2=A_1\cap A_3=A_2\cap A_3=\{e_1\}.$$

However, in this configuration,

$$A_1\cap A_2\cap A_3=\{e_1\}.$$

Thus, both $\omega^{(0)}$ and $\omega^{(1)}$ satisfy the same prescribed matrix $M$, but they have different three-way overlap structures.

**2. Their induced atom vectors**

For three sets, the exact membership atoms are indexed by

$$\emptyset,\ \{1\},\ \{2\},\ \{3\},\ \{1,2\},\ \{1,3\},\ \{2,3\},\ \{1,2,3\}.$$

We write an atom vector in the coordinate order

$$a=\left(a_{\emptyset},a_{\{1\}},a_{\{2\}},a_{\{3\}},a_{\{1,2\}},a_{\{1,3\}},a_{\{2,3\}},a_{\{1,2,3\}}\right).$$

Recall that $a_B$ counts elements whose exact membership set is $B$.

**First overlap pattern: three separate pairwise overlaps**

For $\omega^{(0)}$,

| Element | Exact membership set |
|---|---|
| $e_1$ | $\{1,2\}$ |
| $e_2$ | $\{1,3\}$ |
| $e_3$ | $\{2,3\}$ |
| $e_4$ | $\emptyset$ |
| $e_5$ | $\emptyset$ |

Therefore,

$$a^{(0)} = (2,0,0,0,1,1,1,0).$$

The interpretation is:

$$a_{\emptyset}^{(0)} = 2,$$

because $e_4, e_5$ belong to no sets, and

$$a_{\{1,2\}}^{(0)} = a_{\{1,3\}}^{(0)} = a_{\{2,3\}}^{(0)} = 1,$$

because $e_1, e_2, e_3$, respectively, occupy the three pair-only atoms. In particular,

$$a_{\{1,2,3\}}^{(0)} = 0.$$

**Second overlap pattern: one triple overlap**

For $\omega^{(1)}$,

| Element | Exact membership set |
|---|---|
| $e_1$ | $\{1,2,3\}$ |
| $e_2$ | $\{1\}$ |
| $e_3$ | $\{2\}$ |
| $e_4$ | $\{3\}$ |
| $e_5$ | $\emptyset$ |

Therefore,

$$a^{(1)} = (1,1,1,1,0,0,0,1).$$

The interpretation is:

$$a_{\{1,2,3\}}^{(1)} = 1,$$

because $e_1$ belongs to all three sets, while

$$a_{\{1\}}^{(1)} = a_{\{2\}}^{(1)} = a_{\{3\}}^{(1)} = 1,$$

because $e_2, e_3, e_4$ belong exclusively to $A_1, A_2, A_3$, respectively. Finally,

$$a_{\emptyset}^{(1)} = 1,$$

because $e_5$ belongs to no set.

The two feasible numerical Venn diagrams are therefore

| | $a_{\emptyset}$ | $a_{\{1\}}$ | $a_{\{2\}}$ | $a_{\{3\}}$ | $a_{\{1,2\}}$ | $a_{\{1,3\}}$ | $a_{\{2,3\}}$ | $a_{\{1,2,3\}}$ |
|---|---|---|---|---|---|---|---|---|
| $a^{(0)}$ | 2 | 0 | 0 | 0 | 1 | 1 | 1 | 0 |
| $a^{(1)}$ | 1 | 1 | 1 | 1 | 0 | 0 | 0 | 1. |

### 3. The feasible atom-vector space $\mathcal{A}_2(N, M)$

For the present example, $N = 5$, so $\mathcal{A}_2(N, M) = \mathcal{A}_2(N, M)|_{N=5} = \mathcal{A}_2(5, M)$.

The feasible atom-vector space is

$$\mathcal{A}_2(5, M)|_{N=5} = \{a^{(0)}, a^{(1)}\}.$$

That is, under the prescribed first- and second-order constraints, there are exactly two possible atom-count configurations.

To see why, write

$$r := a_{\{1,2,3\}}.$$

The three pairwise constraints imply

$$a_{\{1,2\}} + r = 1,$$

$$a_{\{1,3\}} + r = 1,$$

and

$$a_{\{2,3\}} + r = 1.$$

Hence,

$$a_{\{1,2\}} = a_{\{1,3\}} = a_{\{2,3\}} = 1 - r.$$

Because all atom counts must be nonnegative integers,

$$r \in \{0, 1\}.$$

- If $r = 0$, then the pairwise intersections are occupied by three distinct elements, producing $a^{(0)}$.
- If $r = 1$, then one element belongs to all three sets, producing $a^{(1)}$.

Thus $\mathcal{A}_2(5, M)$ is the space of two feasible **Venn-diagram count patterns**, not the space of concrete labelled set systems.

## 4. The configuration classes $\boldsymbol{\Omega(a^{(0)})}$ and $\boldsymbol{\Omega(a^{(1)})}$

The full configuration space splits into two disjoint classes:

$$\Omega_2(U, M) = \Omega(a^{(0)}) \cup \Omega(a^{(1)}),$$

where

$$\Omega(a) = \{(A_1, A_2, A_3) \in \Omega_2(U, M)\text{:the induced atom vector is } a\}.$$

**The class $\boldsymbol{\Omega(a^{(0)})}$**

A configuration in $\Omega(a^{(0)})$ has:

- one element in atom $\{1,2\}$;
- one distinct element in atom $\{1,3\}$;
- one distinct element in atom $\{2,3\}$;
- two elements in atom $\emptyset$.

Let

$$e_{12}, e_{13}, e_{23}$$

denote three distinct elements chosen from $U$, where their subscripts indicate their atoms. Then the sets are uniquely determined by

$$A_1 = \{e_{12}, e_{13}\},$$

$$A_2 = \{e_{12}, e_{23}\},$$

$$A_3 = \{e_{13}, e_{23}\}.$$

The remaining two elements of $U$ lie outside all sets.

For instance, the assignment

$$e_{12} = e_1,\ e_{13} = e_2,\ e_{23} = e_3$$

produces $\omega^{(0)}$. But the assignment

$$e_{12} = e_5,\ e_{13} = e_1,\ e_{23} = e_4$$

produces a different labelled configuration:

$$A_1 = \{e_5, e_1\},\ A_2 = \{e_5, e_4\},\ A_3 = \{e_1, e_4\}.$$

The number of configurations in this class is

$$|\ \Omega(a^{(0)})\ | = \frac{5!}{2!\ 1!\ 1!\ 1!} = 60.$$

Equivalently, choose the elements for the three nonempty pair-atoms in order:

$$5 \cdot 4 \cdot 3 = 60.$$

The two remaining elements occupy the outside atom, whose internal ordering does not matter.

**The class $\Omega(a^{(1)})$**

A configuration in $\Omega(a^{(1)})$ has:

- one element in atom $\{1,2,3\}$;
- one element in atom $\{1\}$;
- one element in atom $\{2\}$;
- one element in atom $\{3\}$;
- one element in atom $\emptyset$.

Let

$$e_{123}, e_1^{\circ}, e_2^{\circ}, e_3^{\circ}, e_{\emptyset}$$

be the five elements of $U$, assigned to those five roles. The superscript $e^{\circ}$ merely emphasizes that $e_i^{\circ}$ belongs **only** to $A_i$. Then

$$A_1 = \{e_{123}, e_1^{\circ}\},$$

$$A_2 = \{e_{123}, e_2^{\circ}\},$$

$$A_3 = \{e_{123}, e_3^{\circ}\}.$$

For example, taking

$$e_{123} = e_1,\ e_1^{\circ} = e_2,\ e_2^{\circ} = e_3,\ e_3^{\circ} = e_4,\ e_{\emptyset} = e_5$$

produces $\omega^{(1)}$.

Every atom in this pattern has size one. Therefore,

$$|\Omega(a^{(1)})| = \frac{5!}{1!\ 1!\ 1!\ 1!\ 1!} = 120.$$

**5. The partition function**

The partition function sums the sizes of the two configuration classes:

$$\begin{aligned} Z_{N,M}^{(2)} = Z_{5,M}^{(2)} &= \sum_{a \in \mathcal{A}_2(5,M)} w(a) \\ &= w(a^{(0)}) + w(a^{(1)}) \\ &= 60 + 120 \\ &= 180. \end{aligned}$$

Thus,

$$Z_{5,M}^{(2)} = |\Omega_2(U, M)| = 180.$$

Its meaning is exactly this:

There are 180 labelled triples $(A_1, A_2, A_3)$ of subsets of $U = \{e_1, \ldots, e_5\}$ that satisfy the prescribed set sizes and pairwise-intersection sizes.

The equality

$$Z_{5,M}^{(2)} = |\Omega(a^{(0)})| + |\Omega(a^{(1)})|$$

holds because every feasible labelled configuration induces exactly one atom vector, namely either $a^{(0)}$ or $a^{(1)}$.

**6. The induced distribution on atom vectors**

The second-order MAO null model chooses uniformly from the labelled configuration space:

$$\mathbb{P}_{U,M}^{(2)}(\omega) = \frac{1}{180}, \ \omega \in \Omega_2(U, M).$$

After one forgets which particular labelled elements occupy the atoms and retains only the numerical Venn diagram, the atom-vector distribution becomes

$$\mathbb{P}_{U,M}^{(2)}(a) = \frac{|\Omega(a)|}{|\Omega_2(U, M)|} = \frac{w(a)}{Z_{5,M}^{(2)}}.$$

Consequently,

$$\mathbb{P}^{(2)}_{U,M}\left(a^{(0)}\right) = \frac{60}{180} = \frac{1}{3},$$

whereas

$$\mathbb{P}^{(2)}_{U,M}\left(a^{(1)}\right) = \frac{120}{180} = \frac{2}{3}.$$

Therefore, although both atom vectors satisfy exactly the same prescribed matrix $M$,

$$\mathbb{P}^{(2)}_{U,M}\left(a^{(0)}\right) \neq \mathbb{P}^{(2)}_{U,M}\left(a^{(1)}\right).$$

The vector $a^{(1)}$ is more probable because it corresponds to more ways to allocate the labelled elements $e_1, \dots, e_5$ to the required membership atoms.

**7. Occupancy statistics in the two cases**

The two atom vectors yield the following occupancy profiles:

| | $X_{=0}$ | $X_{=1}$ | $X_{=2}$ | $X_{=3}$ |
|---|---|---|---|---|
| $a^{(0)}$ | 2 | 0 | 3 | 0 |
| $a^{(1)}$ | 1 | 3 | 0 | 1. |

Thus,

$$X_{\geq 3} = X_{=3} = a_{\{1,2,3\}}.$$

Under the uniform labelled-configuration MAO model,

$$\mathbb{P}^{(2)}_{U,M}(X_{\geq 3} = 0) = \frac{1}{3},$$

and

$$\mathbb{P}^{(2)}_{U,M}(X_{\geq 3} = 1) = \frac{2}{3}.$$

This example displays the complete hierarchy:

$$\text{labelled elements } e_1, \dots, e_5 \longrightarrow \text{labelled configurations } (\mathrm{A}_1, \mathrm{A}_2, \mathrm{A}_3) \longrightarrow \text{atom vector } a \longrightarrow \text{occupancy statistics } X_{=t}, X_{\geq t}.$$

The labels $e_1, \dots, e_5$ distinguish concrete configurations, whereas the atom vector records only how many elements have each exact membership pattern.

**Supplementary Material 2**

**S2. Extended Exact Validation**

This supplement reports extended exact-validation benchmarks for the atomic computation. The validation suite includes symmetric and asymmetric margin matrices, direct labelled-set enumeration for small instances, and higher-dimensional atomic calculations. For every atomic benchmark, complete PMFs, factorial moments recovered from these PMFs, and factorial moments accumulated directly over the feasible atomic region were compared using exact rational arithmetic.

For selected benchmarks, the validation further compares the present second-order MAO norm formulas with the corresponding previously implemented formulas and verifies the associated repeated-vector inequalities. Complete machine-readable outputs and executable scripts are supplied with the Supplementary Software.

**S2.1 Independent labelled-set enumeration checks**

For small instances, all labelled set systems satisfying the prescribed matrix $M$ were enumerated directly. These checks include multiple $N = 10$ benchmarks with $T = 3,\ T = 4$, and $T = 5$, together with the symmetric $N = 9, T = 4$ benchmark reported in the main text.

For each such instance, direct labelled enumeration agrees exactly with atom-weighted enumeration for the partition function,

$$| \Omega(N, M) | = Z_{N,M},$$

the complete PMFs of $X_{=t}$ and $X_{\geq t}$, factorial moments through the requested order, and the corresponding MAO norms. These calculations provide an independent check because direct enumeration does not use atomic vectors or multinomial atom weights.

**S2.2 Higher-order asymmetric benchmark: $\boldsymbol{N = 20,\ T = 8}$**

To test the atomic calculation beyond small symmetric examples, we considered the following asymmetric $N = 20,\ T = 8$ benchmark:

$$M = \begin{pmatrix} 4 & 3 & 2 & 0 & 1 & 2 & 1 & 1 \\ 3 & 9 & 4 & 2 & 1 & 4 & 3 & 2 \\ 2 & 4 & 6 & 2 & 1 & 2 & 1 & 0 \\ 0 & 2 & 2 & 7 & 1 & 4 & 5 & 1 \\ 1 & 1 & 1 & 1 & 2 & 1 & 1 & 1 \\ 2 & 4 & 2 & 4 & 1 & 8 & 5 & 1 \\ 1 & 3 & 1 & 5 & 1 & 5 & 10 & 3 \\ 1 & 2 & 0 & 1 & 1 & 1 & 3 & 6 \end{pmatrix}.$$

The diagonal entries specify the eight set sizes, while the off-diagonal entries specify the prescribed pairwise intersections. This benchmark is substantially asymmetric in both its marginal sizes and pairwise overlaps.

All exact-$t$ and threshold statistics were evaluated for

$$t = 0,1, \ldots ,8,$$

corresponding to the nine possible membership orders in a system of $T = 8$ sets. The comparison was carried out through norm order

$$r = 0,1, \ldots ,5.$$

For every $t$ and $r$, the second-order MAO norms obtained from the present formulas agreed exactly with the independently implemented reference formulas. In particular,

$$\| t^r \|_{8;20,M}^{(2)} = \| t^r \|_{8;20,M}^{(2),\mathrm{ref}},$$

and

$$\|[t,8]^r\|_{8;20,M}^{(2)} = \|[t,8]^r\|_{8;20,M}^{(2),ref}.$$

The agreement holds for every occupancy level $t = 0, \ldots ,8$ and every order $r = 0, \ldots ,5$.

For illustration, selected exact-occupancy norms are

$$\| 1^1 \|_{8;20,M}^{(2)} = \frac{45965}{304406},$$

$$\| 2^2 \|_{8;20,M}^{(2)} = \frac{1399678}{14459285},$$

and

$$\| 3^5 \|_{8;20,M}^{(2)} = \frac{786535}{589938828}.$$

Selected threshold norms are

$$\| [1,8]^1 \|_{8;20,M}^{(2)} = \frac{1491777}{1522030},$$

$$\| [2,8]^3 \|_{8;20,M}^{(2)} = \frac{8026364}{14459285},$$

and

$$\| [4,8]^5 \|_{8;20,M}^{(2)} = \frac{5621}{147484707}.$$

**S2.3 more validation examples**

The formulas were validated using three computational routes.

1. **Atom-weighted full PMF calculation.** Each feasible atom vector is weighted by $w(a)$, and the complete PMFs of $X_{=t}$ and $X_{\geq t}$ are accumulated.
2. **Atomic factorial-moment and MAO norm calculation.** The same feasible atomic region is traversed, but the quantities accumulated are the factorial-moment numerators needed for equations (19) and (22). The resulting moments are compared with moments recovered from the complete PMFs.
3. **Direct labelled-set enumeration.** For sufficiently small $N$, all labelled set systems satisfying the imposed matrix $M$ are enumerated directly. This method does not use atom vectors or multinomial atom weights and therefore independently verifies the partition function, distributions, moments, and norms.

The first two routes use the same atomic feasible region but test distinct computational outputs: full distributional aggregation versus direct factorial-moment accumulation. The third route is combinatorially independent and is used where brute-force enumeration is computationally feasible.

**Table S.21. Benchmark cases used for exact validation**

| Case | $N$ | $T$ | Diagonal entries | Pairwise intersections | Direct enumeration |
|---|---|---|---|---|---|
| Symmetric small case | 9 | 4 | (3,3,3,3) | all equal to 1 | Yes |
| Symmetric larger case | 12 | 4 | (4,4,4,4) | all equal to 1 | No |
| Asymmetric case | 12 | 4 | (5,4,6,5) | (2,2,2,2,2,1) | No |

For the symmetric $N = 9, T = 4$ benchmark,

$$M = \begin{pmatrix} 3 & 1 & 1 & 1 \\ 1 & 3 & 1 & 1 \\ 1 & 1 & 3 & 1 \\ 1 & 1 & 1 & 3 \end{pmatrix},$$

the atomic enumeration contains six feasible atom vectors. Its partition function is

$$Z_{9,M} = 808{,}920.$$

Direct enumeration of labelled set systems gives exactly the same value:

$$|\Omega(9, M)| = 808{,}920.$$

The complete exact-$t$ and at-least-$t$ PMFs agree across all three methods, as do factorial moments and MAO norms through order three.

**Table S.22. Selected exact-$t$ results for the $N = 9, T = 4$ symmetric benchmark**

| Statistic | Mean | Variance | Support |
|---|---|---|---|
| $X_{=0}$ | 216/107 | 2136/11449 | {0, 2, 3} |
| $X_{=1}$ | 312/107 | 15648/11449 | {0, 3, 8} |
| $X_{=2}$ | 336/107 | 10368/11449 | {0, 3, 6} |
| $X_{=3}$ | 96/107 | 1056/11449 | {0, 1} |
| $X_{=4}$ | 3/107 | 312/11449 | {0, 1} |

The exact PMFs in this case are

$$\mathbb{P}(X_{=0} = 0,2,3) = \left(\frac{3}{107}, \frac{96}{107}, \frac{8}{107}\right),$$

$$\mathbb{P}(X_{=1} = 0,3,8) = \left(\frac{8}{107}, \frac{96}{107}, \frac{3}{107}\right),$$

$$\mathbb{P}(X_{=2} = 0,3,6) = \left(\frac{3}{107}, \frac{96}{107}, \frac{8}{107}\right),$$

$$\mathbb{P}(X_{=3} = 0,1) = \left(\frac{11}{107}, \frac{96}{107}\right),$$

and

$$\mathbb{P}(X_{=4} = 0,1) = \left(\frac{104}{107}, \frac{3}{107}\right).$$

The threshold distributions follow from the joint atomic law, not from an assumption of independence among the exact occupancy counts.

**Table S.2.3. Selected at-least-$t$ results for the $N = 9, T = 4$ symmetric benchmark**

| Statistic | Mean | Variance | Support |
|---|---|---|---|
| $X_{\geq 0}$ | 9 | 0 | {9} |
| $X_{\geq 1}$ | 747/107 | 2136/11449 | {6, 7, 9} |
| $X_{\geq 2}$ | 435/107 | 6264/11449 | {1, 4, 6} |

| Statistic | Mean | Variance | Support |
|---|---|---|---|
| $X_{\geq 3}$ | $99/107$ | $792/11449$ | $\{0,1\}$ |
| $X_{\geq 4}$ | $3/107$ | $312/11449$ | $\{0,1\}$ |

Several features are visible immediately. First,

$$X_{\geq 1} = N - X_{=0},$$

so $X_{\geq 1}$ and $X_{=0}$ have identical variances. Second, because the constraints permit at most one element with membership order at least three in this benchmark,

$$X_{\geq 3} \in \{0,1\}.$$

Third,

$$X_{\geq 4} = X_{=4}$$

because $T = 4$, as required by equation (16).

For the symmetric $N = 12, T = 4$ benchmark, the exact computation yields

$$\mathbb{E}[X_{\geq 2}] = 5,\ \mathrm{Var}(X_{\geq 2}) = 1,$$

and

$$\mathbb{E}[X_{\geq 3}] = \frac{1}{2},\ \mathrm{Var}(X_{\geq 3}) = \frac{1}{4}.$$

For the asymmetric $N = 12, T = 4$ benchmark, the model produces a nontrivial higher-order overlap structure despite fixed first- and second-order constraints. In particular,

$$\mathbb{E}[X_{\geq 2}] = \frac{20}{3},\ \mathrm{Var}(X_{\geq 2}) = \frac{17}{36},$$

while

$$\mathbb{E}[X_{\geq 3}] = \frac{49}{24},\ \mathrm{Var}(X_{\geq 3}) = \frac{119}{576}.$$

**Table S.2.4. Summary of validation outcomes**

| Case | PMF versus factorial moments | Exact-$t$ norms | At-least-$t$ norms | Direct enumeration |
|---|---|---|---|---|
| $N = 9, T = 4$, symmetric | Exact agreement for $t = 0, \ldots, 4$ | Exact agreement through order 3 | Exact agreement through order 3 | Exact agreement |

| **Case** | **PMF versus factorial moments** | **Exact-$t$ norms** | **At-least-$t$ norms** | **Direct enumeration** |
|---|---|---|---|---|
| $N = 12, T = 4$, symmetric | Exact agreement for $t = 0, \ldots, 4$ | Exact agreement through order 3 | Exact agreement through order 3 | Not run |
| $N = 12, T = 4$, asymmetric | Exact agreement for $t = 0, \ldots, 4$ | Exact agreement through order 3 | Exact agreement through order 3 | Not run |

All computations also satisfy the deterministic identities (1)-(3) exactly. In particular, the PMF-derived expectations satisfy

$$\sum_{t=0}^{T} \mathbb{E}\,[X_{=t}] = N,$$

$$\sum_{t=0}^{T} t\ \mathbb{E}[X_{=t}] = \sum_{i} M_{ii},$$

and

$$\sum_{t=0}^{T} \binom{t}{2} \mathbb{E}[X_{=t}] = \sum_{i<j} M_{ij}.$$